\documentclass[12pt]{amsart}

\usepackage[latin1]{inputenc}

\usepackage{amsmath,amsfonts} 
\usepackage{amssymb}

\usepackage{pdflscape}

\usepackage{pgf,tikz}
\usetikzlibrary{arrows}
\usetikzlibrary{patterns}

\usepackage{hyperref}

\def\+{\!+\!}
\def\-{\!-\!}
\def\g{\mathfrak{g}}

\def\ds{\displaystyle}
\def\sc{\scriptstyle}
\def\ov{\overline}
\def\a{{\bf a}}
\def\b{{\bf b}}
\def\c{{\bf c}}
\def\v{{\bf v}}
\def\x{{\bf x}}
\def\y{{\bf y}}
\def\z{{\bf z}}
\def\ch{{\rm\,ch\,}}
\def\mod{{\rm\,mod\,}}

\def\sgn{\textrm{sgn}}
\def\max{{\rm max}}

\def\cal{\mathcal}

\newtheorem{Theorem}{Theorem}[section]
\newtheorem{Corollary}[Theorem]{Corollary}
\newtheorem{Proposition}[Theorem]{Proposition}
\newtheorem{Conjecture}[Theorem]{Conjecture}
\newtheorem{Definition}[Theorem]{Definition}

\numberwithin{equation}{section}

\newcommand{\N}{{\mathbb N}}
\newcommand{\R}{{\mathbb R}}
\newcommand{\C}{{\mathbb C}}
\newcommand{\Z}{{\mathbb Z}}

\def\qed{\hfill$\Box$\\ \medskip}

\title{Stretched Newell--Littlewood coefficients}

\author{Ronald C. King}

\thanks{
Mathematical Sciences, University of Southampton, 
Southampton SO17 1BJ, England ({\tt r.c.king@soton.ac.uk})}

\begin{document}

\begin{abstract}
Newell--Littlewood coefficients $n_{\mu,\nu}^{\lambda}$ are the multiplicities occurring in the decomposition of 
products of universal characters of the orthogonal and symplectic groups. They may also be expressed, or even defined 
directly in terms of Littlewood--Richardson coefficients, $c_{\mu,\nu}^{\lambda}$. Both 
sets of coefficients have stretched forms $c_{t\mu,t\nu}^{t\lambda}$ and $n_{t\mu,t\nu}^{t\lambda}$,
where $t\kappa$ is the partition obtained by multiplying each part of the partition $\kappa$ by the 
integer $t$. It is known that $c_{t\mu,t\nu}^{t\lambda}$ is a polynomial in $t$ and here it is shown that
$n_{t\mu,t\nu}^{t\lambda}$ is an Ehrhart quasi-polynomial in $t$ with minimum quasi-period at most $2$. 
The evaluation of $n_{t\mu,t\nu}^{t\lambda}$ is effected both by deriving its generating function and 
by establishing a hive model analogous to that used for the calculation of $c_{t\mu,t\nu}^{t\lambda}$. 
These two approaches lead to a whole battery of conjectures about the nature of the quasi-polynomials
$n_{t\mu,t\nu}^{t\lambda}$. These include both positivity, stability and saturation conjectures that are 
supported by a significant amount of data from a range of examples.
\end{abstract}

\maketitle

\section{Introduction}\label{Sec-intro}

There exist irreducible representations, $V_G^\lambda$, of each of the classical Lie groups, $G$, 
whose characters, $\ch V_G^\lambda$, are specified by their highest weights, $\lambda$, which 
take the form of partitions. 
The decomposition of the tensor product of such irreducible 
representations gives rise to multiplicities, $m_{\mu,\nu}^\lambda(G)$,
that are defined at the level of characters by
\begin{equation}\label{eqn-mG}
   \ch V_G^\mu\ \ch V_G^\nu = \sum_{\lambda} m_{\mu,\nu}^{\lambda}(G) \ \ch V_G^\lambda\,.
\end{equation}

In the case of the general linear group, $GL_r$, these multiplicities are known as Littlewood--Richardson 
coefficients and are denoted here by $c_{\mu,\nu}^\lambda$. 
In the case of the orthogonal and symplectic groups, $SO_{2r+1}$, $Sp_{2r}$ and $SO_{2r}$, with $r$ sufficiently large, 
Newell~\cite{N1,N2} and Littlewood~\cite{L2}, by means of quite different arguments, arrived at identical results
allowing the corresponding tensor product multiplicities, $n_{\mu,\nu}^{\lambda}$, that we refer to as Newell--Littlewood 
coefficients, to be expressed in terms of Littlewood--Richardson coefficients as follows: 
\begin{equation}\label{eqn-nl}
   n_{\mu,\nu}^{\lambda}= \sum_{\alpha,\beta,\gamma}\ c_{\alpha,\beta}^\mu \, c_{\alpha,\gamma}^\nu \, c_{\beta,\gamma}^\lambda\,,
\end{equation}
where the sum is over all partitions $\alpha$, $\beta$ and $\gamma$.

These Newell--Littlewood coefficients have been the subject of more recent interest in the hands of Gao~{\it et al.}~\cite{GOY1,GOY2}.
They took (\ref{eqn-nl}) as their definition and made a systematic study of their properties that included some 
remarks and results on stretched Newell--Littlewood coefficients that stimulated the work presented here. 

For each partition $\lambda$ let its weight, i.e. the sum of its parts, be denoted by $|\lambda|$ and its length, 
i.e. the number of its non-zero parts, by $\ell(\lambda)$. It is well known that
\begin{equation}\label{eqn-c-sym}
	 c_{\mu,\nu}^\lambda=c_{\nu,\mu}^\lambda\quad\mbox{and}\quad
	 c_{\mu,\nu}^\lambda=0~~\mbox{unless $|\lambda|=|\mu|\+|\nu|$ and $\ell(\lambda)\leq\ell(\mu)\+\ell(\nu)$}.
\end{equation}
It then follows from (\ref{eqn-nl}) that
\begin{equation}\label{eqn-n-sym}
	\begin{array}{l}
	n_{\mu,\nu}^\lambda=n_{\nu,\mu}^\lambda=n_{\nu,\lambda}^\mu=n_{\lambda,\nu}^\mu=n_{\lambda,\mu}^\nu=n_{\mu,\lambda}^\nu~~\mbox{and}\cr\cr
	n_{\mu,\nu}^\lambda=0~~\mbox{unless $|\lambda|\+|\mu|\+|\nu|$ is even, $|\lambda|\leq|\mu|\+|\nu|$ and $\ell(\lambda)\leq\ell(\mu)\+\ell(\nu)$}\,.\cr
	\end{array}
\end{equation}

For any $t\in\Z_{\geq0}$ and partition $\lambda=(\lambda_1,\lambda_2,\ldots,\lambda_p)$ let $t\lambda$
denote the partition $(t\lambda_1,t\lambda_2,\ldots,t\lambda_p)$. 
Then $c_{t\mu,t\nu}^{t\lambda}$ and $n_{t\mu,t\nu}^{t\lambda}$ are referred to as stretched
Littlewood--Richardson and Newell--Littlewood coefficients, respectively.
Their generating functions take the form:
\begin{equation}\label{eqn-CN}
   C_{\mu,\nu}^\lambda(w)=\sum_{t=0}^\infty c_{t\mu,t\nu}^{t\lambda}\ w^t \quad\mbox{and}\quad 
	 N_{\mu,\nu}^\lambda(w) =\sum_{t=0}^\infty n_{t\mu,t\nu}^{t\lambda}\ w^t\,. 
\end{equation}
It has been established by Rassart~\cite{R} that $c_{t\mu,t\nu}^{t\lambda}$ is a polynomial in $t$, and 
as will be shown here, it follows from the work of De Loera and McAllister~\cite{DeLMcA} that $n_{t\mu,t\nu}^{t\lambda}$
is a quasi-polynomial in $t$ of minimum quasi-period at most $2$, that is to say 
\begin{equation}\label{eqn-PePo}
   c_{t\mu,t\nu}^{t\lambda}=P(t)\quad\mbox{for all $t$}\quad\mbox{and}\quad
	   n_{t\mu,t\nu}^{t\lambda}=\begin{cases}
		                               P_e(t)&\mbox{for $t$ even};\cr
																	 P_o(t)&\mbox{for $t$ odd},\cr
		                          \end{cases}
\end{equation}
where $P(t)$, $P_e(t)$ and $P_o$ are all polynomials in $t$.
In terms of generating functions this implies that~\cite{St}
\begin{equation}\label{eqn-FG}
   C_{\mu,\nu}^\lambda(w)=\frac{F(w)}{(1-w)^d}\quad\mbox{and}\quad 
	 N_{\mu,\nu}^\lambda(w)=\frac{G(w)}{(1-w)^{d_1}(1-w^2)^{d_2}} 
\end{equation}
with $F(w)$ and $G(w)$ polynomials in $w$ of degrees less than $d$ and $d_1+2d_2$, respectively.

For example, in the Littlewood--Richardson case
\begin{equation}\label{eqn-lr-ex}
   c_{(6,5,3),(6,4,1)}^{(9,7,5,4)}=7 \quad\mbox{and}\quad c_{t(6,5,3),t(6,4,1)}^{t(9,7,5,4)}=(t\+1)(5t^2\+10t\+6)/6
\end{equation}
with
\begin{equation}\label{eqn-lr-gf-ex}
    C_{(6,5,3),(6,4,1)}^{(9,7,5,4)}(w)=\frac{w^2\+3w\+1}{(1\-w)^4}\,.
\end{equation}
On the other hand, in the Newell--Littlewood case
\begin{equation}\label{eqn-nl-ex1}
   n_{(5,3),(4,1)}^{(5,2)}=6 \quad\mbox{and}\quad   
	     n_{t(5,3),t(4,1)}^{t(5,2)}=\begin{cases}
		                               (t\+2)(14t^2\+23t\+12)/24&~~t~\mbox{even};\cr
																	 (t\+1)(14t^2\+37t\+21)/24&~~t~\mbox{odd},\cr
		                           \end{cases}
\end{equation}
with
\begin{equation}\label{eqn-nl-gf-ex1}
   N_{(5,3),(4,1)}^{(5,2)}(w)=\ds\frac{3w^2\+3w\+1}{(1\-w)^3(1\-w^2)},
\end{equation}
while
\begin{equation}\label{eqn-nl-ex2}
    n_{(5,3),(4,1)}^{(4,2)}=0 \quad\mbox{and}\quad
    n_{t(5,3),t(4,1)}^{t(4,2)}=\begin{cases}
		                               (t\+2)(19t^2\+40t\+24)/48&~~t~\mbox{even};\cr 
																	 ~~~0&~~t~\mbox{odd},\cr 
		                           \end{cases}
\end{equation}
with
\begin{equation}\label{eqn-nl-gf-ex2}
    N_{(5,3),(4,1)}^{(4,2)}(w)=\ds\frac{7w^4\+11w^2\+1}{(1\-w^2)^4}\,,
\end{equation}

The two Newell--Littlewood examples illustrate some of the properties we wish to explore in an attempt
to find analogues of the following statements that apply to the Littlewood--Richardson case:
\begin{itemize}
\item[\bf LR(i)] Theorem [Knutson-Tao]~\cite{KnTa}: $c_{\mu,\nu}^\lambda>0 \iff c_{t\mu,t\nu}^{t\lambda}>0$,~~~ for all integers $t>0$;
\item[\bf LR(ii)] Theorem [Fulton]~\cite{KTW}: $c_{\mu,\nu}^\lambda=1 \iff c_{t\mu,t\nu}^{t\lambda}=1$,~~~ for all integers $t>0$;
\item[\bf LR(iii)] Theorem [Rassart]~\cite{R}: $c_{t\mu,t\nu}^{t\lambda}$ is a polynomial in $t$ with rational coefficients;
\item[\bf LR(iv)] Conjecture~\cite{KTT}: All coefficients of the polynomial $c_{t\mu,t\nu}^{t\lambda}$ are non-negative.
\item[\bf LR(v)] Conjecture~\cite{KTT}: $F(w)$ is a polynomial in $w$ all of whose coefficients are non-negative integers.
\end{itemize}

In the Newell--Littlewood case it is helpful, as in the above examples, to distinguish between those stretched
coefficients $n_{t\mu,t\nu}^{t\lambda}$ for which $|\lambda|\+|\mu|\+|\nu|$ is either even or odd.
Then the closest analogues of the above statements that we offer as conjectures about stretched
Newell--Littlewood coefficients are as follows:
\begin{Conjecture}\label{Conj-EO}
If $|\lambda|\+|\mu|\+|\nu|$ is even then
\begin{itemize}
\item[\bf E(i)] $n_{\mu,\nu}^{\lambda}>0 \iff n_{t\mu,t\nu}^{t\lambda}>0$,~~~ for all integers $t>0$;
\item[\bf E(ii)] $n_{\mu,\nu}^{\lambda}=1~~\mbox{and}~~n_{2\mu,2\nu}^{2\lambda}=1 \iff n_{t\mu,t\nu}^{t\lambda}=1$,~~~ for all integers $t>1$;
\item[\bf E(iii)] $n_{t\mu,t\nu}^{t\lambda}$ is a quasi-polynomial in $t$ of minimum quasi-period at most $2$ with rational coefficients; 
\item[\bf E(iv)] All coefficients of the quasi-polynomial $n_{t\mu,t\nu}^{t\lambda}$ are non-negative;
\item[\bf E(v)] $G(w)$ is a polynomial in $w$ all of whose coefficients are non-negative integers.
\end{itemize}
If $|\lambda|\+|\mu|\+|\nu|$ is odd then $n_{t\mu,t\nu}^{t\lambda}=0$ for all odd integers $t$ and
\begin{itemize}
\item[\bf O(i)] $n_{2\mu,2\nu}^{2\lambda}>0 \iff n_{2t\mu,2t\nu}^{2t\lambda}>0$,~~~ for all integers $t>0$;
\item[\bf O(ii)] $n_{2\mu,2\nu}^{2\lambda}=1 \iff n_{2t\mu,2t\nu}^{2t\lambda}=1$,~~~ for all integers $t>0$;
\item[\bf O(iii)] $n_{t\mu,t\nu}^{t\lambda}$ is a polynomial in $t^2$ with rational coefficients; 
\item[\bf O(iv)] All coefficients of the polynomial $n_{t\mu,t\nu}^{t\lambda}$ are non-negative;
\item[\bf O(v)] $G(w)$ is a polynomial in $w^2$ all of whose coefficients are non-negative integers, and $d_1=0$.
\end{itemize}
\end{Conjecture}

Here we shall offer a proof of the quasi-polynomial and polynomial conjectures {\bf E(iii)} and {\bf O(iii)}, 
as well as some evidence in support of the remaining conjectures. 
The key point to recognise is that Newell--Littlewood coefficients are nothing other than Clebsch--Gordon coefficients
that govern the decomposition of products of universal characters of the classical orthogonal and symplectic groups,
or equivalently their corresponding simple Lie algebras. This relationship between coefficients is explained in Section~\ref{Sec-char-NL}.
This allows us to conclude immediately that stretched Newell--Littlewood coefficients are quasi-polynomial in nature
by virtue of a proposition to this effect established by De Loera and McAllister~\cite{DeLMcA}. This applies to 
all stretched Clebsch--Gordon coefficients of classical simple Lie algebras. They also established that the minimum
quasi-period of all such quasi-polynomials is at most $2$. This allows us to prove 
a Corollary~\ref{Cor-EOiii} of their proposition that comprises parts {\bf E(iii)} and {\bf O(iii)} of Conjecture~\ref{Conj-EO}. 

This same relationship between Clebsch--Gordon coefficients and Newell--Littlewood coefficients allows us to 
extract from the work of Kapovich and Millson~\cite{KM}, Belkale and Kumar~\cite{BK}, and Sam~\cite{Sa} on the saturation 
problem for the orthogonal and symplectic groups some conjectures on Newell--Littewood coefficients. These are stated
at the close of Section~\ref{Sec-char-NL}. They serve to establish the validity of part {\bf O(i)} of Conjecture~\ref{Conj-EO},
and are consistent with, but weaker than part {\bf E(i)}. 
 
Section~\ref{Sec-gen-fun} is concerned with a new method of calculating the generating function for
stretched Newell--Littlewood coefficients that is based on the use of known generating functions for universal characacters~\cite{L1,K2,FJK}. 
This is motivated by a desire to test further the validity of the positivity conjectures {\bf E(iv)} and {\bf O(iv)} that 
were first formulated in Conjecture~4.7 of~\cite{DeLMcA} within the stretched Clebsch--Gordon coefficient context. 
The required formula applicable to stretched Newell--Littlewood coefficients is provided in Theorem~\ref{The-GF}
and is thereafter exploited by means of Xin's algorithm~\cite{X} to calculate both $N_{\mu,\nu}^\lambda(w)$ and $n_{t\mu,t\nu}^{t\lambda}$, first
for all the examples considered Gao~{\it et al.}~\cite{GOY1}, and in subsequent sections for many other examples.

In order to study the nature of the quasi-polynomials $n_{t\mu,t\nu}^{t\lambda}$ for various triples $(\mu,\nu,\lambda)$, 
in particular their degrees, $\deg n_{t\mu,t\nu}^{t\lambda}$, it is helpful to construct a combinatorial model for the evaluation of Newell--Littlewood
coefficients $n_{\mu,\nu}^{\lambda}$ and to examine the impact of scaling the parts of $\mu$, $\nu$ and $\lambda$ by $t$.
Such a model was used by Gao~{\it et al.}~\cite{GOY1} for their evaluation of Clebsch--Gordon coefficients. It expressed the required
coefficients as the number of integer points in certain BZ-polytopes that had been defined for each of the classical simple
Lie algebras by Berenstein and Zelevinsky~\cite{BZ1,BZ2}. These $BZ$-polytopes are defined by a set of linear 
inequalities and equalities involving the parts of $\mu$, $\nu$ and $\lambda$. Scaling all these parts simultaneously 
by the stretching parameter $t$ allows one to identify stretched Clebsch--Gordon coefficients 
with the Erhart quasi-polynomial~\cite{E,St} of the $BZ$-polytope. 

In the case of the general linear algebra, the $BZ$-polytope is equivalent to the more simply defined hive polytope of 
Knutson and Tao~\cite{KnTa}. This is introduced here in Section~\ref{Sec-hives} as a hive model for the evaluation of 
Littlewood--Richardson coefficients. Then, rather than using the $BZ$-polytopes for the orthogonal and symplectic
algebras to evaluate Newell--Littlewood coefficients, we introduce from first principles a new hive model for their evaluation.
The corresponding polytopes, which we refer to as $K$-polytopes, are again defined by a set of linear 
inequalities and equalities involving the parts of $\mu$, $\nu$ and $\lambda$. Scaling by the stretching parameter $t$ 
in the usual way then allows one to identify stretched Newell--Littlewood coefficients with the Erhart quasi-polynomial
of these $K$-polytopes. This hive model enables one to see that 
\begin{equation}\label{eqn-deg-bound}
   \deg n_{t\mu,t\nu}^{t\lambda} \leq 3n(n-1)/2 \quad\mbox{where $n=\max\{\ell(\mu),\ell(\nu),\ell(\lambda)\}$}\,.
\end{equation}
The use of skeletal graphs~\cite{KTT} in obtaining lower upper bounds for particular triple $(\mu,\nu,\lambda)$ is 
illustrated in three examples, the first of which sheds some light on parts {\bf E(ii)} and {\bf O(ii)} of
Conjecture~\ref{Conj-EO} regarding the conditions under which $n_{t\mu,t\nu}^{t\lambda}=1$.

These bounds are satisfied by explicit calculations of $n_{t\mu,t\nu}^{t\lambda}$ in the following sections
which all, more importantly, support the positivity conjectures of Conjecture~\ref{Conj-EO} as well as giving rise 
to certain new stability conjectures. These arise in Section~\ref{Sec-examples} which commences
with an example in which the quasi-polynomials $n_{t\mu,t\nu}^{t\lambda}$ are evaluated for fixed $\mu=\nu$ but arbitrary $\lambda$.
One stability phenomenon takes the form of the independence of $n_{t\mu,t\nu}^{t\lambda}$ on $a$ in the case
$\mu=(a,\sigma)$, $\nu=(a,\tau)$ and $\lambda=(a,\rho)$ for $a$ sufficiently large. 
Section~\ref{Sec-cubes} involves a further exploration of both the positivity and stability conjectures in the case of Newell--Littlewood 
cubes for which $\mu=\nu=\lambda$, with some further results relegated to Appendix~\ref{Sec-appendix}. 

The concluding Section~\ref{Sec-conclusion} includes some remarks about the connections between the various conjectures.

\section{Universal characters}\label{Sec-char-NL}

The classical Lie groups, $G$, of interest here are the general linear groups, $GL_r$, the odd orthogonal group, $SO_{2r+1}$,
the symplectic group, $Sp_{2r}$, and the even orthogonal group, $SO_{2r}$, for all $r\in\N$. Let ${\cal P}_r$ be the set of all
partitions $\lambda=(\lambda_1,\lambda_2,\ldots,\lambda_r)$ with $\lambda_k\in\Z_{\geq0}$ for $k=1,2,\ldots,r$ and 
$\lambda_k\geq\lambda_{k+1}$ for $k=1,2,\ldots,r-1$. The partition $\lambda$ is said to have length $\ell(\lambda)=p$ if $\lambda_k>0$ 
and $\lambda_k=0$ for all $k>p$. In such a case, one often drops trailing zeros and writes $\lambda=(\lambda_1,\lambda_2,\ldots,\lambda_p)$. 
Each of the groups $G$ possesses a finite-dimensional irreducible representation, $V_G^\lambda$, of highest weight $\lambda$ for 
each $\lambda\in{\cal P}_r$. 

The characters may be evaluated through the use of Weyl's character formula:
\begin{equation}\label{eqn-Weyl}
  \ch V_G^\lambda = \sum_{w\in W_\g} \sgn(w) e^{w(\lambda+\rho)}\ \bigg/\  \prod_{\alpha\in \Delta^+} (e^{\alpha/2}-e^{-\alpha/2})\,,
\end{equation}
where $W_\g$ is the Weyl group of the Lie algebra $\g$ corresponding to the group $G$, $\sgn(w)=\pm1$ is the
signature or parity of $w$, $\Delta_\g^+$ is the set of positive roots of $\g$ and $\rho_\g$ is half the sum of 
the positive roots. It follows that
\begin{equation}\label{eqn-Weyl-exp}
  \prod_{\alpha\in \Delta_\g^+} (1-e^{-\alpha}) \ch V_G^\lambda = \sum_{w\in W_\g} \sgn(w) e^{w(\lambda+\rho_\g)-\rho_\g}\,.  
	                    = e^\lambda + \sum_{\kappa} \left(\pm e^\kappa\right)\,,
\end{equation}
where all terms in the sum on the right are distinct. Each $\kappa=w(\lambda\+\rho_\g)\-\rho_\g$ with $w\neq id$, the identity element of $W_\g$,
does not, unlike $\lambda$, lie in the fundamental, positive Weyl chamber.
Hence, for any coefficients, $m_\mu$
\begin{equation}\label{eqn-weyl-proj}
  [e^\lambda]\ \prod_{\alpha\in \Delta_\g^+} (1-e^{-\alpha})\, \sum_\mu  m_\mu\,\ch V_G^\mu = m_\lambda
\end{equation}
where $[e^\lambda](\cdots)$ is the coefficient of $e^\lambda$ in the expansion of $(\cdots)$.
This allows one to evaluate the coefficients appearing in (\ref{eqn-mG}) by means of the formula
\begin{equation}\label{eqn-cg}
  m_{\mu,\nu}^\lambda(G)= [e^\lambda]\ \prod_{\alpha\in \Delta_\g^+} (1-e^{-\alpha}) \,\ch V_G^\mu \,\ch V_G^\nu\,.
\end{equation}

In each case of interest here the character $\ch V_{G}^\lambda$ may be expressed in terms of a sequence of indeterminates $\x=(x_1,x_2,\ldots)$ whose 
non-vanishing components may be identified either with eigenvalues $x_i$ of group elements or with formal exponentials 
$x_i=e^{\epsilon_i}$ in the Euclidean basis of the root space of the algebra $\g$. For each of our reductive Lie groups, $G$,
the corresponding simple classical Lie algebra, $\g$, and their positive roots are given by:
\begin{equation}\label{eqn-roots}
\begin{array}{|l|l|l|}
\hline
G&\g&\Delta_\g^+\cr
\hline
GL_r&A_{r-1}\dot+D_1&\{\epsilon_i-\epsilon_j\,|\,1\leq i<j\leq r\}\cr
SO_{2r+1}&B_r&\{\epsilon_i\pm\epsilon_j\,|\,1\leq i<j\leq r\}\cup\{\epsilon_i\,|\,1\leq i \leq r\}\cr
Sp_{2r}&C_r&\{\epsilon_i\pm\epsilon_j\,|\,1\leq i<j\leq r\}\cup\{2\epsilon_i\,|\,1\leq i \leq r\}\cr
SO_{2r}&D_r&\{\epsilon_i\pm\epsilon_j\,|\,1\leq i<j\leq r\}\cr
\hline
\end{array}
\end{equation}
It is convenient to set $x_i=e^{\epsilon_i}$ and $\ov{x}_i=x_i^{-1}=e^{-\epsilon_i}$ for $i=1,2,\ldots,r$.
In this Euclidean basis $\lambda=\lambda_1\epsilon_1\+\lambda_2\epsilon_2\+\cdots\+\lambda_r\epsilon_r$ 
and $e^\lambda=\x^\lambda=x_1^{\lambda_1}x_2^{\lambda_2}\cdots x_r^{\lambda_r}$ for any partition $\lambda$ of 
length $\ell(\lambda)\leq r$, with the inclusion where necessary of trailing zeros.
The products in (\ref{eqn-weyl-proj}) take the form
\begin{equation}\label{eqn-tab-delta}
\begin{array}{|l|l|}
\hline
G&\prod_{\alpha\in \Delta_\g^+} (1-e^{-\alpha})\cr
\hline
GL_r&\prod_{1\leq i<j\leq r} (1-\ov{x}_ix_j) \cr
SO_{2r+1}&\prod_{1\leq i<j\leq r} (1-\ov{x}_ix_j) (1-\ov{x}_i\ov{x}_j)\ \prod_{1\leq i\leq r} (1-\ov{x}_i) \cr
Sp_{2r}&\prod_{1\leq i<j\leq r} (1-\ov{x}_ix_j) (1-\ov{x}_i\ov{x}_j) \prod_{1\leq i\leq r} (1-\ov{x}_i^2) \cr
SO_{2r}&\prod_{1\leq i<j\leq r} (1-\ov{x}_ix_j) (1-\ov{x}_i\ov{x}_j)\cr
\hline
\end{array}
\end{equation}

The connection between complex Lie algebras and compact Lie groups is such that the parameters $x_i$ and $\ov{x}_i$, together with $1$, 
may also be interpreted as eigenvalues of group elements, in which case the relevant characters can be specified in the following way:
\begin{equation}\label{eqn-char-spec}
\begin{array}{rcl}
\ch V_{GL_r}^\lambda&=&gl_\lambda(\x) \quad\mbox{with}\quad \x=(x_1,x_2,\ldots,x_r,0,0,\ldots)\cr
\ch V_{SO_{2r+1}}^\lambda&=&oo_\lambda(\x) \quad\mbox{with}\quad \x=(x_1,x_2,\ldots,x_r,\ov{x}_1,\ov{x}_2,\ldots,\ov{x}_r,1,0,0,\ldots) \cr
\ch V_{Sp_{2r}}^\lambda&=&sp_\lambda(\x) \quad\mbox{with}\quad \x=(x_1,x_2,\ldots,x_r,\ov{x}_1,\ov{x}_2,\ldots,\ov{x}_r,0,0,\ldots) \cr
\ch V_{SO_{2r}}^\lambda&=&eo_\lambda(\x) \quad\mbox{with}\quad \x=(x_1,x_2,\ldots,x_r,\ov{x}_1,\ov{x}_2,\ldots,\ov{x}_r,0,0,\ldots) \cr
\end{array}
\end{equation}
Here each character is expressed as a specialisation of 
an appropriate universal character $\g_\lambda(\x)$ with $\x=(x_1,x_2,\ldots)$ a countably infinite sequence.
The corresponding universal characters themselves, without any restrictions on $\x$ are 
defined by means of the generating functions~\cite{L1,K2,FJK}:
\begin{equation}\label{eqn-char-gf}
\begin{array}{rcl}
 \ds \prod_{i,a} (1-x_iy_a)^{-1} &=&\ds \sum_\lambda gl_\lambda(\x)\ gl_\lambda(\y)\,;\cr
 \ds \prod_{i,a} (1-x_iy_a)^{-1} \prod_{a\leq b} (1-y_ay_b) &=&\ds \sum_\lambda oo_\lambda(\x)\ gl_\lambda(\y)\,;\cr
 \ds \prod_{i,a} (1-x_iy_a)^{-1} \prod_{a<b} (1-y_ay_b) &=&\ds \sum_\lambda sp_\lambda(\x)\ gl_\lambda(\y)\,;\cr
 \ds \prod_{i,a} (1-x_iy_a)^{-1} \prod_{a\leq b} (1-y_ay_b) &=&\ds \sum_\lambda eo_\lambda(\x)\ gl_\lambda(\y)\,,\cr
\end{array}
\end{equation}
where $\x=(x_1,x_2,\ldots)$ and $\y=(y_1,y_2,\ldots)$ and the sum is over all partitions $\lambda$.
 
The universal character $gl_\lambda(\x)$ is nothing other than the Schur function $s_\lambda(\x)$
whose product rule~\cite{LR} defines the Littlewood--Richardson coefficients $c_{\mu,\nu}^\lambda$.
The product rules for the universal characters take the form~\cite{N2,L2,K2,FJK}
\begin{equation}\label{eqn-prod-cn}
\begin{array}{rcl}
      gl_\mu(\x)\ gl_\nu(\x) &=&\ds \sum_\lambda c_{\mu,\nu}^\lambda \ gl_\lambda(\x)\,; \cr
			oo_\mu(\x)\ oo_\nu(\x) &=&\ds \sum_\lambda n_{\mu,\nu}^\lambda \ oo_\lambda(\x)\,; \cr
			sp_\mu(\x)\ sp_\nu(\x) &=&\ds \sum_\lambda n_{\mu,\nu}^\lambda \ sp_\lambda(\x)\,; \cr
			eo_\mu(\x)\ eo_\nu(\x) &=&\ds \sum_\lambda n_{\mu,\nu}^\lambda \ eo_\lambda(\x)\,, \cr
\end{array}
\end{equation}
where the coefficients $n_{\mu,\nu}^\lambda$ are precisely the Newell--Littlewood coefficients given by (\ref{eqn-nl}).

The product rules are universal in that they apply in the case of all countably infinite sequences
$\x=(x_1,x_2,\ldots)$. However, under the rank-dependent specialisations of (\ref{eqn-char-spec}) the characters
$\g_\lambda(\x)$ appearing on the right in (\ref{eqn-prod-cn}) may not all be independent. In such a case those for 
which $\ell(\lambda)>r$ are subject to modification rules that can be expressed in a variety of different
ways~\cite{N2,K1,KoTe,BKW}. However, if $r\geq\ell(\mu)\+\ell(\nu)$ then no
modifications are required so that the product rules (\ref{eqn-prod-cn}) apply directly to products 
of characters of $GL_r$, $SO_{2r+1}$, $Sp_{2r}$ and $SO_{2r}$ for all $r\geq\ell(\mu)\+\ell(\nu)$.
As we shall see in Section~\ref{Sec-gen-fun} this observation gives us a way not only of evaluating  
the quasi-polynomial stretched Newell--Littlewood coefficients but also of evaluating their generating functions.

The coefficients $m_{\mu,\nu}^\lambda(G)$ appearing (\ref{eqn-mG}) are referred to variously
as Clebsch--Gordon coefficients, Kronecker product coefficients or tensor product coefficients 
in those cases for which $G$ is a Lie group with corresponding Lie algebra $\g$. To 
distinguish various cases of interest we adopt the following notation:
\begin{equation}\label{eqn-Gg}
\begin{array}{|c|c|l|}
\hline
G&\g&m_{\mu,\nu}^\lambda(G)\cr
\hline
GL_r&A_{r-1}\dotplus D_1&m_{\mu,\nu}^\lambda(gl_r)\cr
SO_{2r+1}&B_r&m_{\mu,\nu}^\lambda(oo_r)\cr
Sp_{2r}&C_r&m_{\mu,\nu}^\lambda(sp_r)\cr
SO_{2r}&D_r&m_{\mu,\nu}^\lambda(eo_r)\ \cr
\hline
\end{array}
\end{equation}
where $r$ is the rank of the relevant finite-dimensional complex simple Lie algebra.
 
Following the observation made previously, 
\begin{equation}\label{eqn-cg-nl}
  m_{\mu,\nu}^\lambda(oo_r)=m_{\mu,\nu}^\lambda(sp_r)=m_{\mu,\nu}^\lambda(eo_r)=n_{\mu,\nu}^\lambda 
\end{equation}
for all $r\geq\ell(\mu)+\ell(\nu)$. 
That is to say, as the rank $r$ increases, each of these Clebsch--Gordon coefficients attains a stable limit 
and this limit coincides with the corresponding Newell--Littlewood coefficient specified by the same 
triple of partition labels $(\mu,\nu,\lambda)$. This result extends to the stretched case:
\begin{equation}\label{eqn-cgt-nlt}
  m_{t\mu,t\nu}^{t\lambda}(oo_r)=m_{t\mu,t\nu}^{t\lambda}(sp_r)=m_{t\mu,t\nu}^{t\lambda}(eo_r)=n_{t\mu,t\nu}^{t\lambda}
\end{equation}
for all $t\in\Z_{\geq0}$ and all $r\geq\ell(\mu)+\ell(\nu)$. 

The behaviour with respect to rank $r$ of stretched Clebsch--Gordon coefficients for each of the Lie algebras
$B_r$, $C_r$ and $D_r$ is illustrated below in the case $\mu=\nu=\lambda=(2,1,1)$ for which it is required 
that $r\geq 3$.

\begin{equation}\label{eqn-cgt-table}
\begin{array}{|c|l|}
\hline
\mathfrak{g}
&\begin{array}{l}
            G(w)/(1\!-\!w)^{d_1}(1\!-\!w^2)^{d_2}\cr
						\begin{cases} P_e(t) & t~\mbox{even}\cr
                          P_o(t) & t~\mbox{odd}\cr
						\end{cases} \cr
 \end{array}\cr
\hline\hline
\begin{array}{c}
B_r\simeq so(2r\+1)\cr r=3\cr
\end{array}
&\begin{array}{l}
                (w^2 \+ w \+ 1)/(1\!-\!w)^3(1\!-\!w^2)\cr
						 \begin{cases}    (t\+2)(2t^2\+5t\+4)/8 & t~\mbox{even}\cr
							 (t\+1)(2t^2+7t\+7)/8 & t~\mbox{odd}\cr
						\end{cases}\cr
\end{array}\cr
\hline
\begin{array}{c}
B_r\simeq so(2r\+1)\cr r\geq 4\cr
\end{array}
&\begin{array}{l}
                 (w^6 \+ w^5 \+ 8w^4 \+ 4w^3 \+ 5w^2 \+ w \+ 1)/(1\!-\!w)^3(1\!-\!w^2)^4\cr
          \begin{cases}      (t \+ 2)^2(t \+ 4)(7t^3 \+ 43t^2 \+ 126t \+ 240)/3840 &t~\mbox{even}\cr
								(t \+ 1)(t \+ 3)(7t^4 \+ 71t^3 \+ 305t^2 \+ 697t \+ 840)/3840 & t~\mbox{odd}\cr
					\end{cases}\cr
\end{array}\cr
\hline
\hline
\begin{array}{c}
C_r\simeq sp(2r)\cr r=3\cr
\end{array}
&\begin{array}{l}
                (w^4 \+ w^2 \+ 1)/(1\!-\!w)(1\!-\!w^2)^3\cr
						 \begin{cases}    (t\+2)(t^2\+4t\+8)/16 & t~\mbox{even}\cr
							 (t\+1)(t^2+2t\+5)/16 & t~\mbox{odd}\cr
						\end{cases}\cr
\end{array}\cr
\hline
\begin{array}{c}
C_r\simeq sp(2r)\cr r\geq 4\cr
\end{array}
&\begin{array}{l}
                 (w^6 \+ w^5 \+ 8w^4 \+ 4w^3 \+ 5w^2 \+ w \+ 1)/(1\!-\!w)^3(1\!-\!w^2)^4\cr
          \begin{cases}      (t \+ 2)^2(t \+ 4)(7t^3 \+ 43t^2 \+ 126t \+ 240)/3840 &t~\mbox{even}\cr
								(t \+ 1)(t \+ 3)(7t^4 \+ 71t^3 \+ 305t^2 \+ 697t \+ 840)/3840 & t~\mbox{odd}\cr
					\end{cases}\cr
\end{array}\cr
\hline
\hline
\begin{array}{c}
D_r\simeq so(2r)\cr r=3\cr
\end{array}
&\begin{array}{l}
                1/(1\!-\!w)\cr
						 \begin{cases} 1& t~\mbox{even}\cr
													 1& t~\mbox{odd}\cr
						\end{cases}\cr
\end{array}\cr
\hline
\begin{array}{c}
D_r\simeq so(2r)\cr r=4\cr
\end{array}
&\begin{array}{l}
                (w^5+ 3w^4 \+ 5w^3 \+ 5w^2 \+ 3w \+ 1)/(1\!-\!w)^4(1\!-\!w^2)^2\cr
						 \begin{cases} (t \+ 2)(6t^4 \+ 33t^3 \+ 104t^2 \+ 152t \+ 80)/160  & t~\mbox{even}\cr
							 (t \+ 1)(6t^4 \+ 39t^3 \+ 131t^2 \+ 229t \+ 155)/160 & t~\mbox{odd}\cr
						\end{cases}\cr
\end{array}\cr
\hline
\begin{array}{c}
D_r\simeq so(2r)\cr r\geq 5\cr
\end{array}
&\begin{array}{l}
                 (w^6 \+ w^5 \+ 8w^4 \+ 4w^3 \+ 5w^2 \+ w \+ 1)/(1\!-\!w)^3(1\!-\!w^2)^4\cr
          \begin{cases}      (t \+ 2)^2(t \+ 4)(7t^3 \+ 43t^2 \+ 126t \+ 240)/3840 &t~\mbox{even}\cr
								(t \+ 1)(t \+ 3)(7t^4 \+ 71t^3 \+ 305t^2 \+ 697t \+ 840)/3840 & t~\mbox{odd}\cr
					\end{cases}\cr
\end{array}\cr
\hline
\hline
\end{array}
\end{equation}

The data underlying the above tabulation were compiled both by using the software package SCHUR~\cite{Sc}
to evaluate the decomposition of products of characters in accordance with (\ref{eqn-mG}) and by using
the well-known character formulae appearing for example in~\cite{FH} and exploiting (\ref{eqn-cg}). 
In this way Clebsch--Gordon coefficients were explicitly calculated for various stretching parameters $t$ for 
each of the Lie algebras of rank $r$, with $r$ ranging from its minimum possible value $3$ to the value $6$
where stability is known to set in. 
The expressions displayed in (\ref{eqn-cgt-table}) were then obtained by fitting the data to 
quasi-polynomials of quasi-period $2$.
The results illustrate the fact that the stable limits of the stretched Clebsch--Gordon 
coefficients of all three Lie algebras coincide, as claimed in (\ref{eqn-cgt-nlt}). 
The fact, that these coefficients differ for lower values of the rank $r$ is a consequence of the modification rules 
applying to universal characters being different for each of the three families of Lie algebras~~\cite{N2,K1,KoTe,BKW}.

This data fitting approach is justified by the fact that De Loera and 
McAllister~\cite{DeLMcA} have already established in their Proposition~1.2
the quasi-polynomial nature of stretched Clebsch--Gordon 
coefficients for all classical Lie algebras, and shown that the minimum quasi-period is at most $2$.
By taking $r$ sufficiently large, that is to say $r\geq\ell(\mu)+\ell(\nu)$ as in (\ref{eqn-cgt-nlt}),
their Proposition~1.2 as applied to Newell--Littlewood coefficients takes the form:
\begin{Proposition}\label{Prop-quasipol}
For all triples of partitions $(\lambda,\mu,\nu)$ and all $t\in\Z_{\geq0}$ the stretched 
Newell--Littlewood coefficients $n_{t\mu,t\nu}^{t\lambda}$ are quasi-polynomials in $t$
with minimum quasi-period at most $2$. 
\end{Proposition}

As a consequence of this we have
\begin{Corollary}\label{Cor-EOiii} \
\begin{itemize}
\item{{\bf E(iii)}} If $|\lambda|\+|\mu|\+|\nu|$ is even and $n_{\mu,\nu}^{\lambda}>0$ then 
$n_{t\mu,t\nu}^{t\lambda}$ is a quasi-polynomial in $t$ with minimum quasi-period at most $2$ with rational coefficients. 
\item{{\bf O(iii)}} If $|\lambda|\+|\mu|\+|\nu|$ is odd and $n_{2\mu,2\nu}^{2\lambda}>0$ then 
$n_{t\mu,t\nu}^{t\lambda}$ is a polynomial in $t^2$ with rational coefficients. 
\end{itemize}
\end{Corollary}

\noindent{\bf Proof}:
As pointed out in the Introduction, Proposition~\ref{Prop-quasipol} is equivalent to the statement that the generating function of (\ref{eqn-CN}) for Newell--Littlewood coefficients takes the form
\begin{equation}\label{eqn-NGw}
\begin{array}{rcl}
    N_{\mu,\nu}^\lambda(w)&=&\ds \sum_{t=0}^\infty n_{t\mu,t\nu}^{t\lambda}\,w^t=\ds\frac{G(w)}{(1\-w)^{d_1}\,(1\-w^2)^{d_2}} \cr 
		&=&\ds\sum_{t\equiv0 (\mod 2)}\,P_e(t)\, w^t\ +\sum_{t\equiv1(\mod 2)}\,P_o(t)\, w^t\,.\cr
\end{array}
\end{equation}
where $G(w)$, $P_e(t)$ and $P_o(t)$ are all polynomials. The fact that $n_{t\mu,t\nu}^{t\lambda}$ is a non-negative integer for all $t$ ensures that the coefficients in $G(w)$, $P_e(t)$ and $P_o(t)$ are rational. 

If $|\lambda|\+|\mu|\+|\nu|$ is even this completes the proof of {\bf E(iii)}. If $|\lambda|\+|\mu|\+|\nu|$ is odd
then $n_{t\mu,t\nu}^{t\lambda}=0$ for all odd $t$, that is to say $P_o(t)=0$. It follows that $d_1=0$
and $G(w)$ must be a polynomial in $w^2$, thereby completing the proof of {\bf O(iii)}.
\qed 

The fact that Newell--Littlewood coefficients are special cases of Clebsch--Gordon coefficients or tensor product multiplicities
of characters of complex reductive Lie groups or their associated Lie algebras also allows one to make considerable progress towards 
the verification of the saturation conjectures {\bf E(i)} and {\bf O(i)}. In particular, Theorem~1.1 of Kapovich and Millson in~\cite{KM},
Theorems~6 and~7 of Belkale and Kumar~\cite{BK}, and Theorem~1.1 of Sam in~\cite{Sa} reveal, when applied to our universal 
characters of the orthogonal and symplectic groups the validity of the following:
\begin{Proposition}\label{Prop-sat}
For all partitions $\mu$, $\nu$ and $\lambda$, and integer $t$ 
\begin{equation}\label{eqn-sat}
    n_{t\mu,t\nu}^{t\lambda}>0 \quad\mbox{for some $t\geq1$} \implies \begin{cases} n_{4\mu,4\nu}^{4\lambda}>0,&~~\mbox{\,\cite{KM}}, \cr
		                                                                      n_{2\mu,2\nu}^{2\lambda}>0,&~~\mbox{\,\cite{BK,Sa}}. \cr 
																																										\end{cases}
\end{equation}
\end{Proposition} 
To justify this one notes that by taking the rank $r$ of the Lie algebra $\g$ of the relevant group $G$
to be such that $r\geq\ell(\mu)+\ell(\nu)\geq \ell(\lambda)$ then $n_{\mu,\nu}^\lambda=m_{\mu,\nu}^\lambda(\g)$,
as explained previously. Moreover, by also taking $r>\ell(\lambda)$ the irreducible representation
of highest weight $\lambda$ will be self-contragredient so that $m_{\mu,\nu}^\lambda(\g)=m_{\mu,\nu,\lambda}^{0}(\g)$, 
where the latter is the multiplicity of the identity representation in the product of the three irreducible representations
of highest weights $\mu$, $\nu$ and $\lambda$. This allows a direct connection to be made with the three-fold case 
of the tensor products encountered in~\cite{KM,BK,Sa}. That the result in~\cite{KM} may be applied to Newell--Littlewood
coefficients follows from the fact that for any three dominant integral weights $\mu$, $\nu$ and $\lambda$,
their sum lies in the root lattice of $B_r$, even though that may not be the case for $C_r$ or $D_r$.

As a consequence of this we have
\begin{Corollary}\label{Cor-EOi} \ 
\begin{itemize}
\item{{\bf E'(i)}} For $|\lambda|\+|\mu|\+|\nu|$ even, if $n_{t\mu,t\nu}^{t\lambda}>0$ for some $t\geq1$, then $n_{2\mu,2\nu}^{2\lambda}>0$: 
\item{{\bf O'(i)}} For $|\lambda|\+|\mu|\+|\nu|$ odd, if $n_{t\mu,t\nu}^{t\lambda}>0$ for some $t\geq1$, then $n_{2\mu,2\nu}^{2\lambda}>0$. 
\end{itemize}
and for $|\lambda|\+|\mu|\+|\nu|$ odd, we have
\begin{itemize}
\item{{\bf O(i)}} $n_{2\mu,2\nu}^{2\lambda}>0 \iff n_{2t\mu,2t\nu}^{2t\lambda}>0$,~~~ for all integers $t>0$ 
\end{itemize}
\end{Corollary}

\noindent{\bf Proof}:
Both {\bf E'(i)} and {\bf O'(i)} are immediate consequences of Proposition~\ref{Prop-sat}. To derive {\bf O(i)}
from {\bf O'(i)} one first notes that for $|\lambda|\+|\mu|\+|\nu|$ odd, $n_{t\mu,t\nu}^{t\lambda}$ may only be $>0$
if $t$ is even. This implies the right to left implication of {\bf O(i)}. The left to right implication 
follows from the fact that $n_{2\mu,2\nu}^{2\lambda}>0$ implies that there exist $\alpha$, $\beta$ and $\gamma$ such that 
$c_{\alpha,\beta}^{2\mu} \, c_{\alpha,\gamma}^{2\nu} \, c_{\beta,\gamma}^{2\lambda}>0$. This implies in turn,
by virtue of {\bf LR(i)}, that $c_{t\alpha,t\beta}^{2t\mu} \, c_{t\alpha,t\gamma}^{2t\nu} \, c_{t\beta,t\gamma}^{2t\lambda}>0$.
\qed

This still leaves open the conjecture {\bf E(i)} that was first proposed by Gao {\it et al.}~\cite{GOY1}. To test this and all remaining 
unproven parts of Conjecture~\ref{Conj-EO} we proceed by gathering together more data.

\section{Generating function approach}\label{Sec-gen-fun}

Although the data of the previous section also supports the positivity conjectures {\bf E(iv)} and {\bf O(iv)} 
on the coefficients of $P_e(t)$ and $P_o(t)$, as well as the positivity conjectures {\bf E(v)} and {\bf O(v)} 
on the coefficients of $G(w)$, these are far from being proved as yet. In fact the positivity of the coefficients
of the polynomials $P_e(t)$ and $P_o(t)$ were first proposed by De Loera and McAllister in Conjecture~4.7 
of~\cite{DeLMcA} within the context of their study of stretched Clebsch--Gordon coefficients. 
This conjecture was based on an extraordinary wealth of data accumulated by the random selection
of hundreds of triples $(\lambda,\mu,\nu)$ involving partitions of impressively high weight. However, the
calculations were only feasible for Lie algebras of comparatively low rank. As a result their tabulation of results 
in support of the positivity conjecture in the Clebsch--Gordon context does not include any Newell--Littlewood
examples. It therefore seems worthwhile calculating some more Newell--Littlewood 
quasi-polynomials from first principles. This may be done by means of the following 
\begin{Theorem}\label{The-GF}
Let $\lambda$, $\mu$ and $\nu$ be partitions of lengths $m$, $p$ and $q$, respectively, with $m\leq p+q$.
Let $\x=(x_1,x_2,\ldots,x_r)$, $\y=(y_1,y_2,\ldots,y_p)$ and $\z=(z_1,z_2,\ldots,z_q)$
with $r=p+q$, and let $\ov{x}_i=x_i^{-1}$, $\ov{y}_i=y_i^{-1}$ and $\ov{z}_i=z_i^{-1}$ for all $i$.
Then
\begin{equation}\label{eqn-n}
        n_{\mu,\nu}^{\lambda}= [\x^{\lambda}\y^{\mu}\z^{\nu}]\ K(\x,\y)K(\x,\z)A(\y)A(\z)C(\ov{\x})V(\x)V(\y)V(\z)
\end{equation}
and
\begin{equation}\label{eqn-Nt}
        N_{\mu,\nu}^\lambda(w)= [\x^{0}\y^{0}\z^{0}]\ \frac{K(\x,\y)K(\x,\z)A(\y)A(\z)C(\ov{\x})V(\x)V(\y)V(\z)}{(1-w/(\x^\lambda\y^\mu\z^\nu))}
\end{equation}
where
\begin{equation}\label{eqn-KCAV}
\begin{array}{l}
\ds K(\x,\y)=\!\!\!\prod_{(i,a)=(1,1)}^{(r,p)}\!\!\!(1-x_iy_a)^{-1}(1-\ov{x}_iy_a)^{-1};~~K(\x,\z)=\!\!\!\prod_{(i,c)=(1,1)}^{(r,p)}\!\!\!(1-x_iz_c)^{-1}(1-\ov{x}_iz_c)^{-1};   \cr\cr
\ds C(\ov{\x})=\prod_{1\leq i\leq j\leq r}(1-\ov{x}_i\ov{x}_j);~~A(\y)=\prod_{1\leq a<b\leq p}(1-y_ay_b);~~A(\z)=\prod_{1\leq c<d\leq q}(1-z_cz_d); \cr\cr
\ds V(\x)=\prod_{1\leq i<j\leq r}(1-\ov{x}_i{x_j});~~V(\y)=\prod_{1\leq a<b\leq p}(1-\ov{y}_ay_b);~~V(\z)=\prod_{1\leq c<d\leq q}(1-\ov{z}_cz_d);\cr
\end{array}
\end{equation}
and $[\x^{\lambda}\,\y^{\mu}\,\z^{\mu}] (\cdots)$ means the coefficient of 
$x_1^{\lambda_1}\cdots x_s^{\lambda_r}y_1^{\mu_1}\cdots y_p^{\mu_p}z_1^{\nu_1}\cdots z_q^{\nu_q}$ in $(\cdots)$.
\end{Theorem}

\noindent{\bf Proof}:
Let $\x=(x_1,x_2,\ldots,x_r,\ov{x}_1,\ov{x}_2,\ldots,\ov{x}_r,0,0,\ldots)$, $\y=(y_1,y_2,\ldots,y_p,0,0,\ldots)$ 
and $\z=(z_1,z_2,\ldots,z_q,0,0,\ldots)$ for the moment. Then the symplectic case of (\ref{eqn-char-gf}) yields,
with the notation of (\ref{eqn-KCAV}),
\begin{equation}\label{eqn-sp2n}
 K(\x,\y)\,A(\y)=\ds \sum_\sigma sp_\sigma(\x)\ gl_\sigma(\y)
\quad\mbox{and}\quad
 K(\x,\z)\,A(\z)=\ds \sum_\tau sp_\tau(\x)\ gl_\tau(\z)\,.
\end{equation}
The choice of $p=\ell(\mu)$ and $q=\ell(\nu)$ in defining $\y$ and $\z$ are the lowest possible values
so as to ensure that the sums over $\sigma$ and $\tau$ include the cases $\sigma=\mu$ and $\tau=\nu$
by virtue of the non-vanishing of $gl_{\mu}(\y)=s_{\mu}(\y)$ and $gl_{\nu}(\z)=s_{\nu}(\z)$. This choice
also ensures that $\ell(\sigma)\leq p$ and $\ell(\tau)\leq q$ in the case of all non-vanishing terms.
Taking the product of these two expressions yields
\begin{equation}\label{eqn-sprho}
K(\x,\y)\,K(\x,\z)\,A(\y)\,A(\z) =
   = \sum_{\rho,\sigma,\tau}  n_{\sigma,\tau}^\rho\ sp_\rho(\x) \ gl_\sigma(\y) \ gl_\tau(\z)\,.
\end{equation}
This time the choice of $r=p+q\geq\ell(\sigma)+\ell(\tau)$ in defining $\x$ ensures 
by virtue of (\ref{eqn-n-sym}) and (\ref{eqn-prod-cn}) that all possible terms $sp_\rho(\x)$ appearing in 
the product of the universal characters $sp_\sigma(\x)$ and $sp_\tau(\x)$ also appear in (\ref{eqn-sprho}).
This includes $sp_{\lambda}(\x)$ since $r=p+q>m$, the length of $\lambda$. 

To project out the coefficients $n_{\mu,\nu}^{\lambda}$ it is necessary to pick 
out of this expression the coefficient of $sp_{\lambda}(\x)\,gl_{\mu}(\y)\,gl_{\nu}(\z)$.
This can be done through the use of (\ref{eqn-weyl-proj}) and the products taken from (\ref{eqn-tab-delta})
for $Sp_{2r}$, $GL_p$ and $GL_q$. 
Again in the notation of (\ref{eqn-KCAV}), these take the form $C(\ov{\x})\,V(\x)$, $V(\y)$ and $V(\z)$.
Putting these results together and taking the coefficient of $[\x^\lambda\,\y^\mu\,\z^\nu]$
yields $n_{\mu,\nu}^{\lambda}$ as in (\ref{eqn-n}), where it is important to note that the
reversion to $\x=(x_1,x_2,\ldots,x_r)$, $\y=(y_1,y_2,\ldots,y_p)$ and $\z=(z_1,z_2,\ldots,z_q)$
is allowed since $\lambda$, $\mu$ and $\nu$ are partitions of lengths $m\leq p+q=r$, $p$ and $q$.

To evaluate $n_{t\mu,t\nu}^{t\lambda}$ for any stretching parameter $t$, one merely has to replace 
$\lambda$, $\mu$ and $\nu$ in (\ref{eqn-n}) by the stretched partitions $t\lambda$, $t\mu$ and $t\nu$, respectively, 
The fact that there is no change in the lengths of these partitions obviates the
necessity of making any changes to the parameters $r$, $p$ and $q$ in (\ref{eqn-KCAV}).

The generating function $N_{\mu,\nu}^{\lambda}(w)$ of (\ref{eqn-Nt}) then follows immediately by expanding $1/(1-w/Z)$
in the form $1+w/Z+w^2/Z^2+\cdots+w^t/Z^t+\cdots$ with $Z^t=\x^{t\lambda}\,\y^{t\mu}\,\z^{t\nu}$.
\qed
 
Similar, but distinct formulae of the same type may be obtained by exploiting the orthogonal cases
of (\ref{eqn-char-gf}) and the $SO(2r)$ and $SO(2r+1)$ cases of (\ref{eqn-tab-delta}). However,
in each case the result is slightly less efficient as a means of calculating $n_{\mu,\nu}^{\lambda}$
and $N_{\mu,\nu}^\lambda(w)$. By exploiting the self-contragredient nature of the relevant irreducible
representations it is also a straightforward task to write down expressions for
$n_{\mu,\nu}^{\lambda}$ and $N_{\mu,\nu}^\lambda(w)$ analogous to those in Theorem~\ref{The-GF}
that are manifestly symmetric in $\mu$, $\nu$ and $\lambda$, namely
\begin{equation}\label{eqn-n-v0}
        n_{\mu,\nu}^{\lambda}= [\v^{0}\x^{\lambda}\y^{\mu}\z^{\nu}]\ 
				K(\v,\x)K(\v,\y)K(\v,\z)A(\x)A(\y)A(\z)C(\ov{\v})V(\v)V(\x)V(\y)V(\z)
\end{equation}
and
\begin{equation}\label{eqn-Nt-v0}
        N_{\mu,\nu}^\lambda(w)= [\v^{0}\x^{0}\y^{0}\z^{0}]\ 
				\frac{K(\v,\x)K(\v,\y)K(\v,\z)A(\x)A(\y)A(\z)C(\ov{\v})V(\v)V(\x)V(\y)V(\z)}{(1-w/(\v^0\x^\lambda\y^\mu\z^\nu))}
\end{equation}
where $\v=(v_1,v_2,\ldots,v_s)$ with $s=p+q+r$. 
However, the use of these expressions is again computationally more demanding than the use of those in Theorem~\ref{The-GF}.  

To exploit Theorem~\ref{The-GF} it is convenient to use Xin's algorithm~\cite{X} as 
implemented in the Maple package Ell2.mpl. This proceeds by successively picking out the constant
terms in each $x_i$, each $y_a$ and each $z_c$. To do this it is necessary to multiply $x_i$ and $\ov{x}_i$ 
in $K(\x,\y)$ by a parameter $u$, and in $K(\x,\z)$ by a parameter $v$, only setting $u=v=1$ after the
elimination of $x_i$. The result is an explicit expression for $N_{\mu,\nu}^\lambda(w)$ in the
form $G(w)/(1-w)^{d_1}(1-w^2)^{d_2}$ from which can be extracted polynomial expressions $P_e(t)$
and $P_o(t)$ for $n_{t\mu,t\nu}^{t\lambda}$ in the cases $t$ even and odd, respectively. 

The following results of Gao~{\it et al.} in section 5.4 of~\cite{GOY1} have been confirmed by this means.
\begin{equation}\label{eqn-GOY}
\begin{array}{|c|l|}
\hline
\begin{array}{c}
\mu\cr\nu\cr\lambda\cr
\end{array}
&\begin{array}{l}
            G(w)/(1\!-\!w)^{d_1}(1\!-\!w^2)^{d_2}\cr
						\begin{cases} P_e(t) & t~\mbox{even}\cr
                          P_o(t) & t~\mbox{odd}\cr
						\end{cases} \cr
 \end{array}\cr
\hline
\begin{array}{c}
(1,1)\cr(1,1)\cr(1,1)\cr
\end{array}
&\begin{array}{l}
                1/(1\!-\!w)(1\!-\!w^2)\cr
						 \begin{cases}    (t\+2)/2 & t~\mbox{even}\cr
							 (t\+1)/2 & t~\mbox{odd}\cr
						\end{cases}\cr
\end{array}\cr
\hline
\begin{array}{c}
(2,1,1)\cr(2,1,1)\cr(1,1,1,1)\cr
\end{array}
&\begin{array}{l}   
                 1/(1\!-\!w)^2(1\!-\!w^2)^2\cr
							 \begin{cases}(t \+ 2)(t \+ 3)(t \+4)/24 & t~\mbox{even}\cr
								(t \+ 1)(t \+ 3)(t \+5)/24 & t~\mbox{odd}\cr
				 \end{cases}\cr			   \end{array}\cr
\hline
\begin{array}{c}
(2,1,1)\cr(2,1,1)\cr(2,1,1)\cr
\end{array}
&\begin{array}{l}
                 (w^6 \+ w^5 \+ 8w^4 \+ 4w^3 \+ 5w^2 \+ w \+ 1)/(1\!-\!w)^3(1\!-\!w^2)^4\cr
          \begin{cases}      (t \+ 2)^2(t \+ 4)(7t^3 \+ 43t^2 \+ 126t \+ 240)/3840 &t~\mbox{even}\cr
								(t \+ 1)(t \+ 3)(7t^4 \+ 71t^3 \+ 305t^2 \+ 697t \+ 840)/3840 & t~\mbox{odd}\cr
					\end{cases}\cr
\end{array}\cr
\hline
\end{array}
\end{equation}

\section{Hive model}\label{Sec-hives}

In the case of the Lie algebra $gl(n)$ the Littlewood--Richardson coefficients may be evaluated using the hive model
introduced by Knutson and Tao~\cite{KnTa}, with properties described in more detail by Buch~\cite{B}.
An integer $n$-hive is a labelling of the vertices of a planar, equilateral triangular graph of side length $n$ 
with integers $a_{ij}$, for $0\leq i\leq j\leq n$,
satisfying certain rhombus inequalities which are to be applied to each elementary rhombus formed 
from the union of any pair of elementary triangles having a common edge whatever their orientation.
One such hive is illustrated below on the left in the case $n=4$.
The elementary triangles and rhombi are shown on the right. They are provided with both vertex and edge labels~\cite{KTT,TKA}.
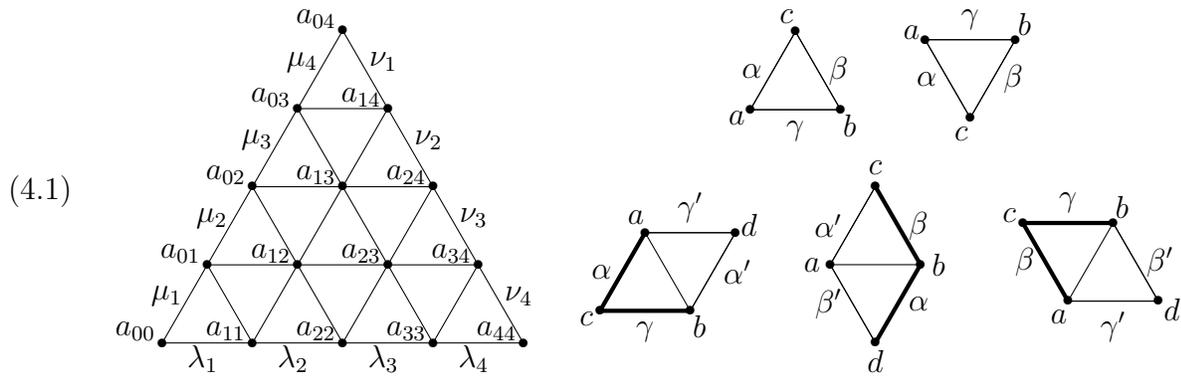
\begin{equation}\label{eqn-LRhive}
\vcenter{\hbox{\begin{tikzpicture}[x={(1cm*0.6,-1.732 cm*0.6)},y={(1cm*0.6,1.732 cm*0.6)}] 
\foreach \i in{0,...,3} \draw(\i,\i)--(\i,4); 
\foreach \j in{1,...,4} \draw(0,\j)--(\j,\j);  
\foreach \i in{0,...,4} \draw(0,\i)--(4-\i,4); 
\foreach\j in{0,...,4} \draw(0-0.37,\j-0.22)node(a0\j){$a_{0\j}$};
\foreach\j in{1,...,4} \draw(1-0.37,\j-0.22)node(a1\j){$a_{1\j}$};
\foreach\j in{2,...,4} \draw(2-0.37,\j-0.22)node(a2\j){$a_{2\j}$};
\foreach\j in{3,...,4} \draw(3-0.37,\j-0.22)node(a3\j){$a_{3\j}$};
\foreach\j in{4,...,4} \draw(4-0.37,\j-0.22)node(a4\j){$a_{4\j}$};
\foreach\j in{0,...,4} \draw(0,\j)node(b0\j){$\sc\bullet$};
\foreach\j in{1,...,4} \draw(1,\j)node(b1\j){$\sc\bullet$};
\foreach\j in{2,...,4} \draw(2,\j)node(b2\j){$\sc\bullet$};
\foreach\j in{3,...,4} \draw(3,\j)node(b3\j){$\sc\bullet$};
\foreach\j in{4,...,4} \draw(4,\j)node(b4\j){$\sc\bullet$};
\foreach\j in{1,...,4} \draw(0-0.25,\j-0.65)node(mu\j){$\mu_{\j}$};
\foreach\j in{1,...,4} \draw(\j-0.35,4+0.25)node(nu\j){$\nu_{\j}$};
\foreach\j in{1,...,4} \draw(\j-0.45,\j-0.65)node(la\j){$\lambda_{\j}$};
\end{tikzpicture}}}
\begin{array}{cc}  
\begin{array}{cc}  
\vcenter{\hbox{\begin{tikzpicture}[x={(1cm*0.6,-1.732 cm*0.6)},y={(1cm*0.6,1.732 cm*0.6)}]
\draw(0,0)node(a){$\sc\bullet$}; \draw(1,1)node(b){$\sc\bullet$}; \draw(0,1)node(c){$\sc\bullet$}; 
\draw(0,1)--(1,1);\draw(0,1)--(0,0);\draw(1,1)--(0,0);
\draw(0-0.1,0-0.2)node{$a$};\draw(1+0.2,1)node{$b$};\draw(0-0.2,1)node{$c$};
\foreach \p/\q/\anchor/\label in{{(0,0)}/{(0,1)}/left/{\alpha},{(0,1)}/{(1,1)}/right/{\beta},{(0,0)}/{(1,1)}/below/{\gamma}} \draw[-]\p--node[\anchor]{$\label$}\q;
\end{tikzpicture}}}
&
\vcenter{\hbox{\begin{tikzpicture}[x={(1cm*0.6,-1.732 cm*0.6)},y={(1cm*0.6,1.732 cm*0.6)}]
\draw(1,0)node(a){$\sc\bullet$}; \draw(1,1)node(b){$\sc\bullet$}; \draw(0,0)node(c){$\sc\bullet$}; 
 \draw(1,0)--(1,1);\draw(1,0)--(0,0);\draw(1,1)--(0,0);
\draw(1,0-0.2)node{$c$};\draw(1,1+0.2)node{$b$};\draw(0-0.2,0-0.1)node{$a$};
\foreach \p/\q/\anchor/\label in{{(0,0)}/{(1,0)}/left/{\alpha},{(1,1)}/{(1,0)}/right/{\beta},{(0,0)}/{(1,1)}/above/{\gamma}} \draw[-]\p--node[\anchor]{$\label$}\q;
\end{tikzpicture}}}
\cr
\end{array} 
\cr
\begin{array}{ccc} 
\vcenter{\hbox{\begin{tikzpicture}[x={(1cm*0.6,-1.732 cm*0.6)},y={(1cm*0.6,1.732 cm*0.6)}]
\draw(0,1)node(a){$\sc\bullet$}; \draw(1,1)node(b){$\sc\bullet$}; \draw(0,0)node(c){$\sc\bullet$}; \draw(1,2)node(d){$\sc\bullet$};
\draw(0,1)--(1,1);\draw[ultra thick](0,1)--(0,0);\draw(0,1)--(1,2);\draw[ultra thick](1,1)--(0,0);\draw(1,1)--(1,2);
\draw(0-0.2,1)node{$a$};\draw(1+0.2,1)node{$b$};\draw(0-0.1,0-0.2)node{$c$};\draw(1+0.1,2+0.2)node{$d$};
\foreach \p/\q/\anchor/\label in{{(0,0)}/{(0,1)}/left/{\alpha},{(1,1)}/{(1,2)}/right/{\alpha'},{(0,0)}/{(1,1)}/below/{\gamma},{(0,1)}/{(1,2)}/above/{\gamma'}} \draw[-]\p--node[\anchor]{$\label$}\q;
\end{tikzpicture}}}
&
\vcenter{\hbox{\begin{tikzpicture}[x={(1cm*0.6,-1.732 cm*0.6)},y={(1cm*0.6,1.732 cm*0.6)}]
 \draw(0,0)node(a){$\sc\bullet$}; \draw(1,1)node(b){$\sc\bullet$}; \draw(0,1)node(c){$\sc\bullet$}; \draw(1,0)node(d){$\sc\bullet$};
 \draw(0,0)--(1,1);\draw(0,0)--(0,1);\draw(0,0)--(1,0);\draw[ultra thick](1,1)--(0,1);\draw[ultra thick](1,0)--(1,1);
\draw(0-0.2,0-0.2)node{$a$};\draw(1+0.2,1+0.2)node{$b$};\draw(0-0.12,1+0.12)node{$c$};\draw(1+0.12,0-0.12)node{$d$};
\foreach \p/\q/\anchor/\label in{{(1,0)}/{(1,1)}/right/{\alpha},{(0,0)}/{(0,1)}/left/{\alpha'},{(0,0)}/{(1,0)}/left/{\beta'},{(0,1)}/{(1,1)}/right/{\beta}} \draw[-]\p--node[\anchor]{$\label$}\q;

\end{tikzpicture}
}}
&
\vcenter{\hbox{\begin{tikzpicture}[x={(1cm*0.6,-1.732 cm*0.6)},y={(1cm*0.6,1.732 cm*0.6)}]
\draw(1,0)node(a){$\sc\bullet$}; \draw(1,1)node(b){$\sc\bullet$}; \draw(0,0)node(c){$\sc\bullet$}; \draw(2,1)node(d){$\sc\bullet$};
 \draw(1,0)--(1,1);\draw[ultra thick](1,0)--(0,0);\draw(1,0)--(2,1);\draw[ultra thick](1,1)--(0,0);\draw(1,1)--(2,1);
\draw(1,0-0.2)node{$a$};\draw(1,1+0.2)node{$b$};\draw(0-0.2,0-0.1)node{$c$};\draw(2+0.2,1+0.1)node{$d$};
\foreach \p/\q/\anchor/\label in{{(0,0)}/{(1,0)}/left/{\beta},{(1,1)}/{(2,1)}/right/{\beta'},{(0,0)}/{(1,1)}/above/{\gamma},{(1,0)}/{(2,1)}/below/{\gamma'}} \draw[-]\p--node[\anchor]{$\label$}\q;
\end{tikzpicture}}}
\cr
\end{array}
\cr
\end{array}
\end{equation}
An edge between vertices labelled $a$ and $b$ is given the label $b-a$ if $b$ is to the right of $a$. 
Thus in each of the two elementary triangles we have 
\begin{equation}\label{eqn-triangle-conditions}
  \alpha=c-a,~~\beta=b-c,~~\gamma=b-a \quad\mbox{so that}\quad \gamma=\alpha+\beta.
\end{equation}
For all three elementary rhombi, the rhombus constraints take the form
\begin{equation}\label{eqn-rhombus-conditions}
   a+b\geq c+d \quad\mbox{so that}\quad \alpha\geq\alpha',~~\beta\geq\beta',~~\gamma\geq\gamma'\,.
\end{equation}
This is indicated in the above diagram by making the edge with the potentially larger edge label thicker than the other 
for each pair of parallel edges.

As emphasised elsewhere~\cite{TKA}, one consequence of the rhombus constraints is
that if they are paired together as in the following three diagrams they imply the
betweenness conditions as specified below each diagram: 
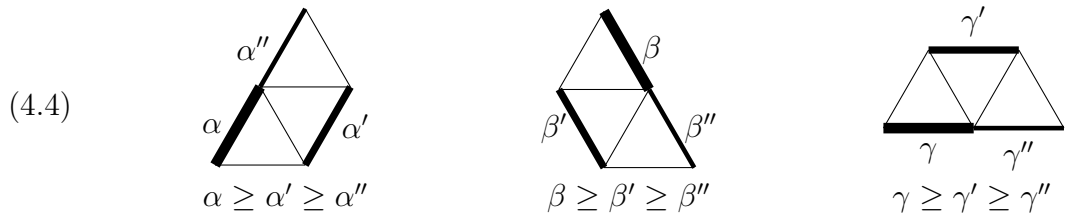
\begin{equation}\label{eqn-betweenness}
\begin{array}{ccc}
\vcenter{\hbox{%
\begin{tikzpicture}
 [x={(1cm*0.6,-1.732 cm*0.6)},y={(1cm*0.6,1.732 cm*0.6)}]
 \foreach \p/\q/\anchor/\label/\pos/\width in {
  {(0,0)}/{(0,1)}/left/\alpha/0.5/4pt,
  {(0,1)}/{(0,2)}/left/\alpha''/0.5/1.75pt,
  {(1,1)}/{(1,2)}/right/\alpha'/0.5/2.75pt
 } \draw[line width=\width]\p--node[pos=\pos,\anchor,black]{$\label$}\q;
 \foreach \p/\q in {
  {(0,1)}/{(1,1)}, 
  {(0,2)}/{(1,2)}, 
  {(0,0)}/{(1,1)}, 
  {(0,1)}/{(1,2)}
 } \draw[-] \p--\q;
\end{tikzpicture}}}
\qquad&\qquad 
\vcenter{\hbox{%
\begin{tikzpicture}
 [x={(1cm*0.6,-1.732 cm*0.6)},y={(1cm*0.6,1.732 cm*0.6)}]
 \foreach \p/\q/\anchor/\label/\pos/\width in{
  {(0,1)}/{(1,1)}/left/\beta' /0.5/2.75pt,
  {(0,2)}/{(1,2)}/right/\beta/0.5/4pt,
  {(1,2)}/{(2,2)}/right/\beta''/0.5/1.75pt
 } \draw[line width=\width]\p--node[pos=\pos,\anchor,black]{$\label$}\q;
 \foreach \p/\q in {
  {(0,1)}/{(0,2)}, 
  {(1,1)}/{(1,2)}, 
  {(1,1)}/{(2,2)}, 
  {(0,1)}/{(1,2)}
 } \draw[-]\p--\q;
\end{tikzpicture}}}
\qquad&\qquad
\vcenter{\hbox{%
\begin{tikzpicture}
 [x={(1cm*0.6,-1.732 cm*0.6)},y={(1cm*0.6,1.732 cm*0.6)}]
 \foreach \p/\q/\anchor/\label/\pos/\width in {
  {(0,0)}/{(1,1)}/below/\gamma/0.5/4pt,
  {(1,1)}/{(2,2)}/below/\gamma''/0.5/1.75pt,
  {(0,1)}/{(1,2)}/above/\gamma'/0.5/2.75pt
 } \draw[line width=\width]\p--node[pos=\pos,\anchor,black]{$\label$}\q;
 \foreach \p/\q in {
  {(0,0)}/{(0,1)}, 
  {(1,1)}/{(1,2)}, 
  {(0,1)}/{(1,1)}, 
  {(1,2)}/{(2,2)}
 } \draw[-]\p--\q;
\end{tikzpicture}}}
\cr
\alpha\geq\alpha'\geq\alpha''\qquad&\qquad\beta\geq\beta'\geq\beta''\qquad&\qquad\gamma\geq\gamma'\geq\gamma''
\cr
\end{array}
\end{equation}
This immediately implies that the sequence of edge labels along any line parallel to one or other of the boundaries
of the hive, including the boundaries themselves, constitute a partition if read from bottom left to top right parallel
to the left-hand boundary, or from top left to bottom right parallel to the right-hand boundary or from left to right 
parallel to the bottom boundary.

\begin{Definition}\label{Def-Hn}
For fixed integer $n$, and partitions $\mu$, $\nu$ and $\lambda$ of lengths 
$\ell(\mu),\ell(\nu),\ell(\lambda)\leq n$, let ${\cal H}^{(n)}(\mu,\nu;\lambda)$ be the set of 
integer $n$-hives $H$, satisfying the constraints (\ref{eqn-triangle-conditions})
and (\ref{eqn-rhombus-conditions}) whose boundary edge labels are specified by the parts of the partitions 
$\mu$, $\nu$ and $\lambda$ in accordance with the formulae
\begin{equation}\label{Eq-LRhive}
  \mu_i=a_{0,i}-a_{0,i-1},~~\nu_i=a_{i,n}-a_{i-1,n},~~\lambda_i=a_{i,i}-a_{i-1,i-1}\quad\mbox{for}\quad i=1,2,\ldots,n\,.
\end{equation}
\end{Definition}
This labelling of boundary edges is illustrated for the case $n=4$ in the above diagram (\ref{eqn-LRhive}). 

With this definition, Fulton has established in the Appendix of~\cite{B} a bijection between such hives and 
the tableaux whose enumeration serves to evaluate the Littlewood--Richardson coefficients $c_{\mu,\nu}^{\lambda}$.
This implies the validity of the following hive model formula for these coefficients as given in Appendix~2 of~\cite{KnTa}.
\begin{Proposition}\label{Prop-LRcoeff-hive}
Let $\mu$, $\nu$ and $\lambda$ be partitions of lengths $\ell(\mu),\ell(\nu),\ell(\lambda)\leq n$.
Then 
\begin{equation}\label{eqn-LRcoeff-hive}
    c_{\mu,\nu}^{\lambda}= \#\{H\in{\cal H}^{(n)}(\mu,\nu;\lambda)\}\,.
\end{equation}
\end{Proposition}

It is then a simple matter to exploit the definition (\ref{eqn-nl}) of Newell--Littlewood coefficients
in terms of Littlewood--Richardson coefficients to arrive at a hive model formula for the former. To this
end, let $K$ be the composite $n$-hive constructed from three standard $n$-hives as shown below:
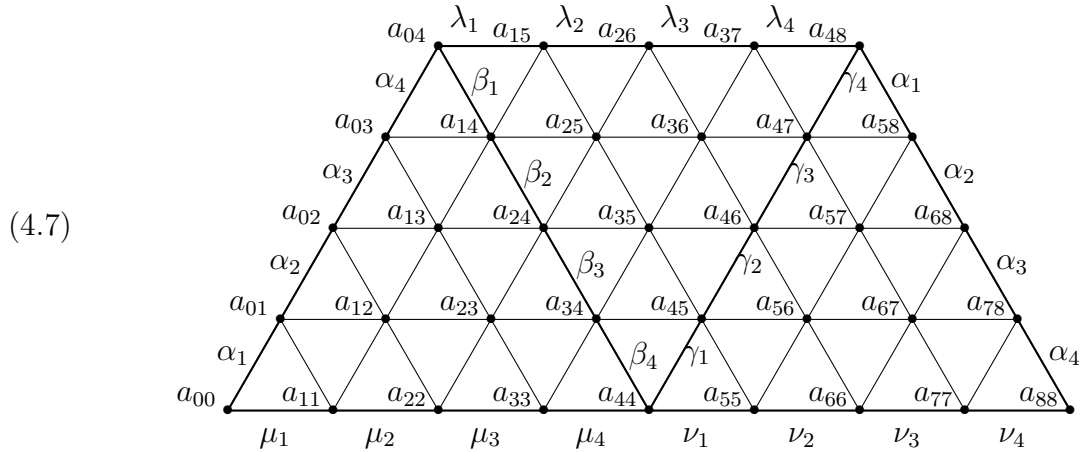
\begin{equation}\label{eqn-NLhive}
\vcenter{\hbox{\begin{tikzpicture}[x={(1cm*0.7,-1.732 cm*0.7)},y={(1cm*0.7,1.732 cm*0.7)}] 
\foreach \i in{0,...,4} \draw(\i,\i)--(\i,4+\i); \foreach \i in{5,...,8} \draw(\i,\i)--(\i,8); \draw[thick](0,0)--(0,4); \draw[thick](4,4)--(4,8);
\foreach \j in{1,...,4} \draw(0,\j)--(\j,\j);  \foreach \j in{5,...,8} \draw(\j-4,\j)--(\j,\j); \draw[thick](0,4)--(4,4); \draw[thick](4,8)--(8,8); 
\foreach \i in{0,...,4} \draw(0,\i)--(8-\i,8); \draw[thick](0,0)--(8,8); \draw[thick](0,4)--(4,8); 
\foreach\j in{0,...,4} \draw(0,\j)node(a0\j){$\sc\bullet$};
\foreach\j in{1,...,5} \draw(1,\j)node(a1\j){$\sc\bullet$};
\foreach\j in{2,...,6} \draw(2,\j)node(a2\j){$\sc\bullet$};
\foreach\j in{3,...,7} \draw(3,\j)node(a3\j){$\sc\bullet$};
\foreach\j in{4,...,8} \draw(4,\j)node(a4\j){$\sc\bullet$};
\foreach\j in{5,...,8} \draw(5,\j)node(a0\j){$\sc\bullet$};
\foreach\j in{6,...,8} \draw(6,\j)node(a1\j){$\sc\bullet$};
\foreach\j in{7,...,8} \draw(7,\j)node(a2\j){$\sc\bullet$};
\foreach\j in{8,...,8} \draw(8,\j)node(a3\j){$\sc\bullet$};
\foreach\j in{0,...,4} \draw(0-0.37,\j-0.22)node(a0\j){$a_{0\j}$};
\foreach\j in{1,...,5} \draw(1-0.37,\j-0.22)node(a1\j){$a_{1\j}$};
\foreach\j in{2,...,6} \draw(2-0.37,\j-0.22)node(a2\j){$a_{2\j}$};
\foreach\j in{3,...,7} \draw(3-0.37,\j-0.22)node(a3\j){$a_{3\j}$};
\foreach\j in{4,...,8} \draw(4-0.37,\j-0.22)node(a4\j){$a_{4\j}$};
\foreach\j in{5,...,8} \draw(5-0.37,\j-0.22)node(a5\j){$a_{5\j}$};
\foreach\j in{6,...,8} \draw(6-0.37,\j-0.22)node(a6\j){$a_{6\j}$};
\foreach\j in{7,...,8} \draw(7-0.37,\j-0.22)node(a7\j){$a_{7\j}$};
\foreach\j in{8,...,8} \draw(8-0.37,\j-0.22)node(a8\j){$a_{8\j}$};
\foreach\j in{1,...,4} \draw(0-0.25,\j-0.65)node(alphaleft\j){$\alpha_{\j}$};
\foreach\j in{1,...,4} \draw(4-0.25+0.4,4+\j-0.65+0.4)node(gamma\j){$\gamma_{\j}$};
\foreach\j in{1,...,4} \draw(4+\j-0.35,8+0.25)node(alpharight\j){$\alpha_{\j}$};
\foreach\j in{1,...,4} \draw(0+\j-0.35,4+0.25)node(beta\j){$\beta_{\j}$};
\foreach\j in{1,...,4} \draw(\j-0.4,\j-0.7)node(mu\j){$\mu_{\j}$};
\foreach\j in{1,...,4} \draw(4+\j-0.4,4+\j-0.7)node(nu\j){$\nu_{\j}$};
\foreach\j in{1,...,4} \draw(\j-0.9,4+\j-0.6)node(la\j){$\lambda_{\j}$};
\end{tikzpicture}}}
\end{equation} 
Such a trapezoidal composite $n$-hive is constructed from three constituent triangular integer $n$-hives, $H_\mu$, $H_\nu$ and $H^\lambda$,
with horizontal boundary edge labels specified by the parts of $\mu$, $\nu$ and $\lambda$, respectively. The $\lambda$ 
triangular $n$-hive has been turned upside down, so that it shares common boundaries with both 
the $\mu$ and $\nu$ triangular $n$-hives. We refer to these internal boundaries as the $\beta$ and $\gamma$ boundaries, respectively,
and the sloping outer boundaries as left and right $\alpha$-boundaries. 

\begin{Definition}\label{Def-K-hive}
For some fixed positive integer $n$, and partitions $\mu$, $\nu$ and $\lambda$ of lengths 
$\ell(\mu)$, $\ell(\nu)$ and $\ell(\lambda)$ all $\leq n$, let ${\cal K}^{(n)}(\mu,\nu;\lambda)$ be the 
set of all composite $n$-hives $K$ with horizontal lower
boundary edge labels specified by the parts of $\mu$ followed by the parts of $\nu$ and horizontal
upper boundary edge labels specified by the parts of $\lambda$, all read from left to right, such that:
(i) the triangle rule (\ref{eqn-triangle-conditions}) is satisfied by each elementary triangle; 
(ii) the rhombus constraints (\ref{eqn-rhombus-conditions}) apply to each elementary rhombus that 
does not cross either the $\beta$- or the $\gamma$-boundary; 
(iii) the bottom to top edge labels on the left $\alpha$-boundary coincide with the top to bottom edge 
labels on the right $\alpha$-boundary.
\end{Definition} 

With this definition we have:
\begin{Proposition}\label{Prop-NLcoeff-hive}
Let $\mu$, $\nu$ and $\lambda$ be partitions of lengths $\ell(\mu),\ell(\nu),\ell(\lambda)\leq n$.
Then 
\begin{equation}\label{eqn-NRcoeff-hive}
    n_{\mu,\nu}^{\lambda}= \#\{K\in{\cal K}^{(n)}(\mu,\nu;\lambda)\}\,.
\end{equation}
\end{Proposition}
 
\noindent{\bf Proof}: First it should be noted that the sets of triangle conditions (\ref{eqn-triangle-conditions}) 
and rhombus constraints (\ref{eqn-rhombus-conditions}) are unchanged by turning the triangles and rhombi 
upside down. In enumerating all possible composite hives $K$ with horizontal boundary edges determined 
by $\mu$, $\nu$ and $\lambda$ the rhombus constraints alone on $H_\mu$, $H_\nu$ and $H^\lambda$ would lead to hives 
with $\alpha$-, $\beta$- and $\gamma$-boundary edge labels specified by all possible partitions $\alpha,\beta,\gamma,\alpha',\beta',\gamma'$
such that $c_{\alpha,\beta'}^{\mu}$, $c_{\gamma,\alpha'}^{\nu}$ and $c_{\beta,\gamma'}^{\lambda}$ are non-zero.
However, the juxtaposition of $H_\mu$, $H_\nu$ and $H^\lambda$ imposes the constraints $\beta'=\beta$ and $\gamma'=\gamma$,
while the above requirement (iii) in the definition of a composite $n$-hive requires that $\alpha'=\alpha$.
The enumeration of such composite $n$-hives therefore yields the product of the Littlewood--Richardson coefficients 
$c_{\alpha,\beta}^{\mu}$, $c_{\gamma,\alpha}^{\nu}$ and $c_{\beta,\gamma}^{\lambda}$, with of course 
$c_{\gamma,\alpha}^{\nu}=c_{\alpha,\gamma}^{\nu}$ for all possible partitions $\alpha$, $\beta$ and $\gamma$
such that these coefficients are non-zero. The formula (\ref{eqn-nl}) for $n_{\mu,\nu}^{\lambda}$
involving a sum over partitions $\alpha$, $\beta$ and $\gamma$ then immediately gives (\ref{eqn-NRcoeff-hive}).
\qed

In this model, the enumeration of all composite hives $K\in{\cal K}^{(n)}(\mu,\nu;\lambda)$ involves
edge labels specified by the parts of $\alpha$, $\beta$ and $\gamma$ that are restricted in
the first instance only by the constraints
\begin{equation}\label{eqn-alpha-beta-gamma}
   2|\alpha|=|\mu|+|\nu|-|\lambda|, \quad 2|\beta|=|\lambda|+|\mu|-|\nu|, \quad 2|\gamma|=|\nu|+|\lambda|-|\mu|\,,
\end{equation}
as can be established by repeated use of the triangle conditions (\ref{eqn-triangle-conditions}), 
but also seen directly of course from (\ref{eqn-nl}).
Since $\alpha$, $\beta$  and $\gamma$ are partitions, these constraints (\ref{eqn-alpha-beta-gamma}) 
suffice to show that $n_{\mu,\nu}^{\lambda}=0$ if $|\mu|+|\nu|+|\lambda|$ is odd. 

\begin{Definition}\label{Def-K-polytope}
For a triple of partitions $(\mu,\nu,\lambda)$ with $\{\ell(\mu),\ell(\nu),\ell(\lambda)\}\leq n$,
let $k=(3n+2)(n+1)/2$. Then $k$ is the number of vertices of the composite $n$-hive $K$ of edge length $n$
with lower and upper edge labels specified by the parts of $(\mu,\nu)$ and $\lambda$, respectively.
Then the $K$-polytope ${\cal P}_{\mu,\nu}^{\lambda}\subset \R^k$ is the convex hull of the points 
$a_{i,j}\in\R^k$ for $j=0,1,\ldots,n$ and $i=0,1,\ldots,2n-j$ 
subject to the linear equalities and inequalities:
\begin{itemize}
  \item{(i)} $a_{0,0}=0,~~\mbox{and}~~a_{0,n}\-a_{0,0}=(|\mu|\+|\nu|\-|\lambda|)/2$;
  \item{(ii)} $a_{n+j,2n}\-a_{n,2n}=a_{0,j}\-a_{0,0}$ for $j=1,2,\ldots,n-1$;
	\item{(iii)} $\mu_i=a_{i,i}\-a_{i-1,i-1}$, $\nu_i=a_{n+i,n+i}\-a_{n+i-1,n+i-1}$, $\lambda_i=a_{i,n+i}\-a_{i-1,n+i-1}$ for $i=1,2,\ldots,n$;
  \item{(iv)} $R_{i,j},U_{i,j},L_{i,j},R_{i+1,j+n},U_{i,j+n-1},L_{i,j+n-1},R_{i+n,j+n},U_{i+n,j+n},L_{i+n,j+n}$ all $\geq0$ for $1\leq i<j\leq n$, where
    $R_{i,j}=a_{i-1,j-1}\+a_{i,j-1}\-a_{i,j}\-a_{i-1,j-2}$, $U_{i,j}=a_{i-1,j-1}\+a_{i,j}\-a_{i-1,j}\-a_{i,j-1}$, $L_{i,j}=a_{i-1,j-1}\+a_{i+1,j}\-a_{i,j}\-a_{i,j-1}$.
\end{itemize}
\end{Definition}
The first part $a_{0,0}=0$ of (i) is just an anchoring condition, with all other conditions involving differences of the form $a_{i,j}-a_{k,l}$.
The second part of (i) is the constraint on $|\alpha|$ given in (\ref{eqn-alpha-beta-gamma}), while (ii) is the left and right boundary matching
condition. Part (iii) corresponds to the lower and upper boundary conditions, while part (iv) corresponds to the rhombus conditions (\ref{eqn-rhombus-conditions})
applied to the three different orientations of a rhombus in each of the three constituent $n$-hives $H_\mu$, $H_\nu$ and $H^\lambda$.

We then have:
\begin{Theorem}\label{The-Ehrhart}
With the notation of Definition~\ref{Def-K-polytope} we have
\begin{equation}\label{eqn-Ehrhart}
   n_{\mu,\nu}^\lambda = i({\cal P}_{\mu,\nu}^{\lambda},1):=\#\{{\cal P}_{\mu,\nu}^{\lambda}\cap \Z^k\} \quad\mbox{and}\quad 
	 n_{t\mu,t\nu}^{t\lambda} = i({\cal P}_{\mu,\nu}^{\lambda},t):=\#\{{\cal P}_{t\mu,t\nu}^{t\lambda}\cap \Z^k\}
\end{equation}
where $i({\cal P},t)$ is the Ehrhart quasi-polynomial of the polytope ${\cal P}$. 
\end{Theorem}
\noindent{\bf Proof}: The first claim is just a restatement of the fact that $n_{\mu,\nu}^\lambda$ 
is the number of composite $K$-hives with vertices labelled by integers $a_{i,j}$ subject to the given
equalities and inequalities applicable to ${\cal P}_{\mu,\nu}^{\lambda}$. All of these can be expressed in the form 
$E(\a)=\b$ and $I(\a)\leq\c$ for some rectangular matrices $E$ and $I$ with integer elements, where $\a$
is a sequence of length $k$ of the labels $a_{i,j}$, and $\b$ and $\c$ are linear in the parts of $\mu$, $\nu$ and $\lambda$.
Scaling these parts by $t$ serves to specify ${\cal P}_{t\mu,t\nu}^{t\lambda}$ whose intersection with $\Z^k$
necessarily yields $n_{t\mu,t\nu}^{t\lambda}$.
\qed 
  
Propositions~\ref{Prop-LRcoeff-hive} and \ref{Prop-NLcoeff-hive} imply that both Littlewood-Richardson and
Newell--Littlewood coefficients are independent of the parameter $n$ that determines the boundary edge lengths of the
$n$-hives and the composite $n$-hives, respectively, provided that these edge lengths are large enough to accommodate 
the numbers of non-vanishing parts of the partitions $\mu$, $\nu$ and $\lambda$. This can be seen within the context of the hive models
themselves. In the case of any $c_{\mu,\nu}^\lambda>0$, if $\lambda$ has a trailing $0$, then so must both $\mu$ and $\nu$,
and the corresponding hive takes the typical form:
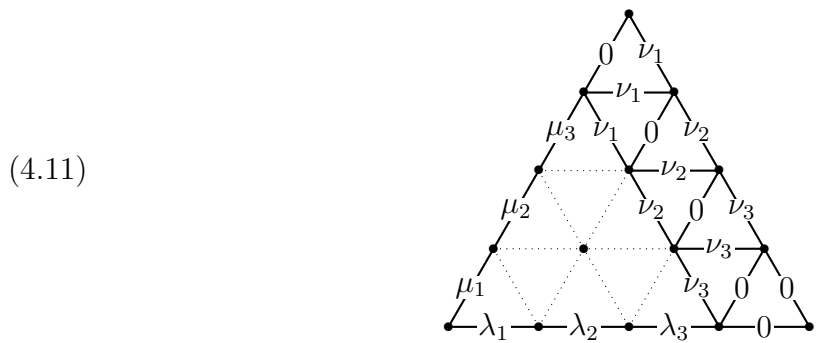
\begin{equation}\label{eqn-Hmn}
\vcenter{\hbox{\begin{tikzpicture}[x={(1cm*0.6,-1.732 cm*0.6)},y={(1cm*0.6,1.732 cm*0.6)}] 
\foreach \i in {1,...,3} \draw[dotted](\i,\i)--(\i,3);
\foreach \i in {1,...,3} \draw[thick](\i,3)--(\i,4);
\draw[thick](0,0)--(0,4); 
\foreach \j in{1,2} \draw[dotted](0,\j)--(\j,\j);  
\foreach \j in{3,4} \draw[thick](0,\j)--(\j,\j); 
\foreach \i in{1,2} \draw[dotted](0,\i)--(3-\i,3);
\foreach \i in{1,...,3} \draw[thick](\i-1,3)--(\i,4);
\draw[thick](0,0)--(4,4);
\foreach\j in{0,...,4} \draw(0,\j)node(a0\j){$\sc\bullet$};
\foreach\j in{1,...,4} \draw(1,\j)node(a1\j){$\sc\bullet$};
\foreach\j in{2,...,4} \draw(2,\j)node(a2\j){$\sc\bullet$};
\foreach\j in{3,...,4} \draw(3,\j)node(a3\j){$\sc\bullet$};
\foreach\j in{4,...,4} \draw(4,\j)node(a4\j){$\sc\bullet$};
\foreach\j in{1,...,3} \draw(\j-0.5,\j-0.5)node[fill=white,inner sep=0.5pt]{$\lambda_\j$};
\foreach\j in{1,...,3} \draw(0,\j-0.5)node[fill=white,inner sep=0.5pt]{$\mu_\j$};
\foreach\j in{1,...,3} \draw(\j-0.5,3)node[fill=white,inner sep=0.5pt]{$\nu_\j$};
\foreach\j in{1,...,3} \draw(\j-0.5,4)node[fill=white,inner sep=0.5pt]{$\nu_\j$};
\foreach\j in{1,...,3} \draw(\j-0.5,4-0.5)node[fill=white,inner sep=0.5pt]{$\nu_\j$};
\foreach\j in{4} \draw(0,\j-0.5)node[fill=white,inner sep=0.5pt]{$0$};
\foreach\j in{4} \draw(1,\j-0.5)node[fill=white,inner sep=0.5pt]{$0$};
\foreach\j in{4} \draw(2,\j-0.5)node[fill=white,inner sep=0.5pt]{$0$};
\foreach\j in{4} \draw(3,\j-0.5)node[fill=white,inner sep=0.5pt]{$0$};
\foreach\j in{4} \draw(\j-0.5,4)node[fill=white,inner sep=0.5pt]{$0$};
\foreach\j in{4} \draw(\j-0.5,\j-0.5)node[fill=white,inner sep=0.5pt]{$0$};
\end{tikzpicture}}}
\end{equation}
The fact that $\mu$ and $\nu$ must also have a trailing $0$ is a consequence of the triangle and rhombus conditions 
(\ref{eqn-triangle-conditions}) and (\ref{eqn-rhombus-conditions}), respectively. These conditions also
suffice to fix all edge labels of the ladder on the right, with its rungs labelled $0$ and its side and remaining interior 
labels specified by the parts of $\nu$. This reduces the enumeration of $4$-hives to that of $3$-hives.
In this way the evaluation of $c_{\mu,\nu}^\lambda$ can always be reduced to the enumeration of all $n$-hives
with $n=\ell(\lambda)$.

The same type of argument applies to the evaluation of $n_{\mu,\nu}^\lambda$, which can be reduced 
to the enumeration of composite $n$-hives with $n=\max\{\ell(\mu),\ell(\nu),\ell(\lambda)\}$. The 
reduction process is illustrated in the following example
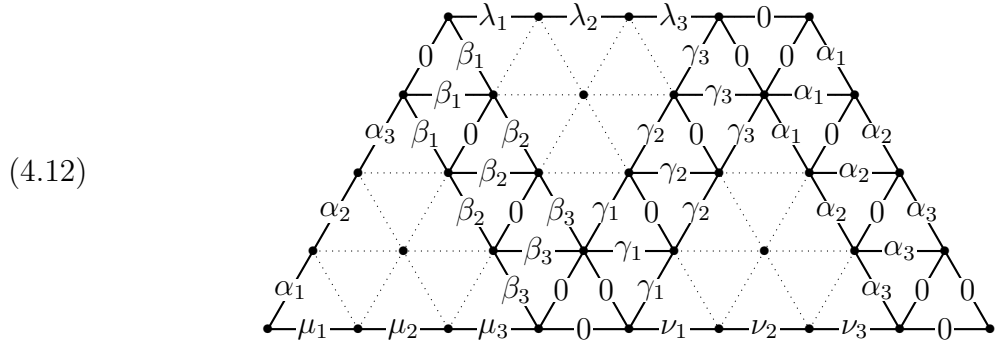
\begin{equation}\label{eqn-K-mn} 
\vcenter{\hbox{\begin{tikzpicture}[x={(1cm*0.6,-1.732 cm*0.6)},y={(1cm*0.6,1.732 cm*0.6)}] 
\foreach \i in {0,...,4} \draw[dotted](\i,\i)--(\i,\i+4); \foreach \i in {4,...,7} \draw[dotted](\i,\i)--(\i,8); 
\draw[thick](0,0)--(0,4); 
\draw[thick](1,3)--(1,4);
\draw[thick](2,3)--(2,4);
\draw[thick](3,3)--(3,4);\draw[thick](3,4)--(3,7);
\draw[thick](4,4)--(4,7);\draw[thick](4,7)--(4,8);
\draw[thick](5,7)--(5,8);
\draw[thick](6,7)--(6,8);
\draw[thick](7,7)--(7,8);
\foreach \j in{1,...,4} \draw[dotted](0,\j)--(\j,\j);  \foreach \j in{4,...,8} \draw[dotted] (\j-4,\j)--(\j,\j); 
\draw[thick](0,3)--(3,3); 
\draw[thick](0,4)--(4,4); 
\draw[thick](3,5)--(4,5); 
\draw[thick](3,6)--(4,6); 
\draw[thick](3,7)--(7,7);
\draw[thick](4,8)--(8,8);
\draw[thick](0,0)--(8,8);
\foreach \i in{1,...,4} \draw[dotted](0,\i)--(8-\i,8);
\foreach \i in{1,...,3} \draw[thick](\i-1,3)--(\i,4);
\foreach \i in{1,...,3} \draw[thick](3,3+\i)--(4,4+\i);
\foreach \i in{1,...,3} \draw[thick](\i+3,7)--(\i+4,8);
\draw[thick](0,4)--(4,8); 
\foreach\j in{0,...,4} \draw(0,\j)node(a0\j){$\sc\bullet$};
\foreach\j in{1,...,5} \draw(1,\j)node(a1\j){$\sc\bullet$};
\foreach\j in{2,...,6} \draw(2,\j)node(a2\j){$\sc\bullet$};
\foreach\j in{3,...,7} \draw(3,\j)node(a3\j){$\sc\bullet$};
\foreach\j in{4,...,8} \draw(4,\j)node(a4\j){$\sc\bullet$};
\foreach\j in{5,...,8} \draw(5,\j)node(a0\j){$\sc\bullet$};
\foreach\j in{6,...,8} \draw(6,\j)node(a1\j){$\sc\bullet$};
\foreach\j in{7,...,8} \draw(7,\j)node(a2\j){$\sc\bullet$};
\foreach\j in{8,...,8} \draw(8,\j)node(a3\j){$\sc\bullet$};
\foreach\j in{1,2,3} \draw(\j-0.5,\j-0.5)node[fill=white,inner sep=0.5pt]{$\mu_\j$};
\foreach\j in{1,2,3} \draw(\j+4-0.5,\j+4-0.5)node[fill=white,inner sep=0.5pt]{$\nu_\j$};
\foreach\j in{1,2,3} \draw(\j-0.5,\j+4-0.5)node[fill=white,inner sep=0.5pt]{$\lambda_\j$};
\foreach\j in{1,...,3} \draw(0,\j-0.5)node[fill=white,inner sep=0.5pt]{$\alpha_\j$};
\foreach\j in{1,...,3} \draw(4+\j-0.5,8)node[fill=white,inner sep=0.5pt]{$\alpha_\j$};
\foreach\j in{1,...,3} \draw(4+\j-0.5,8-0.5)node[fill=white,inner sep=0.5pt]{$\alpha_\j$};
\foreach\j in{1,...,3} \draw(4+\j-0.5,7)node[fill=white,inner sep=0.5pt]{$\alpha_\j$};
\foreach\j in{1,...,3} \draw(\j-0.5,3)node[fill=white,inner sep=0.5pt]{$\beta_\j$};
\foreach\j in{1,...,3} \draw(\j-0.5,4)node[fill=white,inner sep=0.5pt]{$\beta_\j$};
\foreach\j in{1,...,3} \draw(\j-0.5,4-0.5)node[fill=white,inner sep=0.5pt]{$\beta_\j$};
\foreach\j in{1,...,3} \draw(3,4+\j-0.5)node[fill=white,inner sep=0.5pt]{$\gamma_\j$};
\foreach\j in{1,...,3} \draw(4,4+\j-0.5)node[fill=white,inner sep=0.5pt]{$\gamma_\j$};
\foreach\j in{1,...,3} \draw(4-0.5,4+\j-0.5)node[fill=white,inner sep=0.5pt]{$\gamma_\j$};
\foreach\j in{4} \draw(0,\j-0.5)node[fill=white,inner sep=0.5pt]{$0$};
\foreach\j in{4} \draw(1,\j-0.5)node[fill=white,inner sep=0.5pt]{$0$};
\foreach\j in{4} \draw(2,\j-0.5)node[fill=white,inner sep=0.5pt]{$0$};
\foreach\j in{4} \draw(3,\j-0.5)node[fill=white,inner sep=0.5pt]{$0$};
\foreach\j in{8} \draw(4,\j-0.5)node[fill=white,inner sep=0.5pt]{$0$};
\foreach\j in{8} \draw(5,\j-0.5)node[fill=white,inner sep=0.5pt]{$0$};
\foreach\j in{8} \draw(6,\j-0.5)node[fill=white,inner sep=0.5pt]{$0$};
\foreach\j in{8} \draw(7,\j-0.5)node[fill=white,inner sep=0.5pt]{$0$};
\foreach\j in{4} \draw(\j-0.5,4)node[fill=white,inner sep=0.5pt]{$0$};
\foreach\j in{4} \draw(\j-0.5,5)node[fill=white,inner sep=0.5pt]{$0$};
\foreach\j in{4} \draw(\j-0.5,6)node[fill=white,inner sep=0.5pt]{$0$};
\foreach\j in{4} \draw(\j-0.5,7)node[fill=white,inner sep=0.5pt]{$0$};
\foreach\j in{8} \draw(\j-0.5,8)node[fill=white,inner sep=0.5pt]{$0$};
\foreach\j in{4,8} \draw(\j-0.5,\j-0.5)node[fill=white,inner sep=0.5pt]{$0$};
\foreach\j in{8} \draw(\j-4-0.5,\j-0.5)node[fill=white,inner sep=0.5pt]{$0$};
\end{tikzpicture}}}
\end{equation} 
Removal of the three ladders, all of whose edges are fixed by a combination of
triangle and rhombus conditions, reduces the composite $4$-hive to a composite $3$-hive.

Dealing with trailing zeros common to all three partitions $\mu$, $\nu$ and $\lambda$
by the iteration of this process enables any composite $n$-hive to be reduced to one for which 
$n=\max\{\ell(\mu),\ell(\nu),\ell(\lambda)\}$. This immediately implies that
the degree of the corresponding Ehrhart quasi-polynomial $i({\cal P}_{\mu,\nu}^{\lambda},t)$
is bounded by $k$ with $k=(3n+2)(n+1)/2$. This maximum degree can be reduced still further 
by using all the equalities in the definition of the $K$-polytope. This can be done
by noting that, for given $\mu$, $\nu$ and $\lambda$, the parts (i)--(iii) of Definition~\ref{Def-K-polytope} 
fix the vertex labels on the top and bottom horizontal boundary and the right hand diagonal boundary.
It follows that the maximum degree of $n_{t\mu,t\nu}^{t\lambda}$ is $3n(n-1)/2$, that is to say $0,3,9,18,...$ for $n=1,2,3,4,\ldots$.

That this upper bound on the degree of $n_{t\mu,t\nu}^{t\lambda}$ is not reached in any of the three examples of (\ref{eqn-GOY}) is clear.
However, a lower upper bound on the degree may found by exploiting the notion of skeletal graphs~\cite{KTT}. 
In the present context a skeletal graph is a partial labelling
of the underlying graph of the appropriately shaped composite $K$-hive in which edges are identified for which the
labels are fixed from a knowledge of the parts of $\mu$, $\nu$ and $\lambda$, including trailing zeros. These are shown in the examples
below as thick lines along with their labels as determined by the triangle conditions and rhombus constraints. Since the
labels on the right-hand boundary are determined by those on the left-hand boundary, the edges of the
former are shown as dashed lines. All remaining edges are represented as dotted lines to indicate that their labels are not fixed.
The remaining degrees of freedom in enumerating the composite $K$-hives contributing to $n_{t\mu,t\nu}^{t\lambda}$ is
then the number of components of the skeletal graph, including isolated vertices, that are disconnected from the lower and upper 
horizontal boundaries.

All this is illustrated in the case of the three examples of (\ref{eqn-GOY}) as follows:

For $\mu=\nu=\lambda=(1,1)$ the skeletal $K$-hive takes the form:
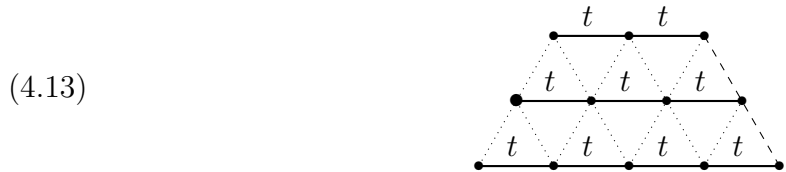
\begin{equation}\label{eqn-Kskel-ex1}
\vcenter{\hbox{\begin{tikzpicture}[x={(1cm*0.5,-1.732 cm*0.5)},y={(1cm*0.5,1.732 cm*0.5)}] 
\draw[dotted](0,0)--(0,2); \draw[dotted](1,1)--(1,3); \draw[dotted](2,2)--(2,4); \draw[dotted](3,3)--(3,4);
\draw[dotted](0,1)--(1,1); \draw[dotted](0,2)--(2,2); \draw[dotted](1,3)--(3,3); \draw[dashed](2,4)--(4,4); 
\draw[thick](0,0)--(4,4); \draw[thick](0,1)--(3,4); \draw[thick](0,2)--(2,4);
\foreach\j in{0,...,2} \draw(0,\j)node(a0\j){$\sc\bullet$};
\foreach\j in{1,...,3} \draw(1,\j)node(a1\j){$\sc\bullet$};
\foreach\j in{2,...,4} \draw(2,\j)node(a2\j){$\sc\bullet$};
\foreach\j in{3,...,4} \draw(3,\j)node(a0\j){$\sc\bullet$};
\foreach\j in{4,...,4} \draw(4,\j)node(a1\j){$\sc\bullet$};
\draw(0,1)node(01){$\bullet$};
\foreach\j in{1,...,2} \draw(\j-0.7,\j-0.4)node(mu\j){$t$};
\foreach\j in{3,...,4} \draw(\j-0.7,\j-0.4)node(nu\j){$t$};
\foreach\j in{1,...,2} \draw(\j-0.7,2+\j-0.4)node(la\j){$t$};
\foreach\j in{1,...,3} \draw(\j-0.7,1+\j-0.4)node(ka\j){$t$};
\end{tikzpicture}}}
\end{equation} 
The skeletal graph contains a single component disconnected from the horizontal boundaries and this is signified by marking 
its left-most vertex by means of an enlarged bullet point. There is therefore just one degree of freedom in enumerating the 
corresponding composite $K$-hives, so the quasi-polynomial $n_{t(1,1),t(1,1)}^{t(1,1)}$ is of degree $1$ in
agreement with (\ref{eqn-GOY}). In this particular case it is easy to go further. The left-hand boundary edge labels
$\alpha=(\alpha_1,\alpha_2)$ must be such that with $\alpha_1+\alpha_2=2t$ with $t\geq\alpha_1\geq\alpha_2$, leading 
to the conclusion that $n_{t(1,1),t(1,1)}^{t(1,1)}=(t+2)/2$ if $t$ is even, and $(t+1)/2$ if $t$ is odd, as 
first found by Gao~{\it et al.}~\cite{GOY1}. 
It can also be seen in this case that $n_{\mu,\nu}^{\lambda}=1$ and $n_{2\mu,2\nu}^{2\lambda}=2$, thereby ruling out any
notion that $n_{\mu,\nu}^{\lambda}=1$ might always imply $n_{t\mu,t\nu}^{t\lambda}=1$ for all $t>0$, as in {\bf LR(ii)}.
However, as yet, we have uncovered no counterexample to the conjecture that, taken together, 
$n_{\mu,\nu}^{\lambda}=1$ and $n_{2\mu,2\nu}^{2\lambda}=1$ imply $n_{t\mu,t\nu}^{t\lambda}=1$ for all $t>0$, as 
would be required by {\bf E(ii)}.
 
In the case $\mu=\nu=(2,1,1)$ and $\lambda=(1,1,1,1)$ the skeletal $K$-hive takes the form:
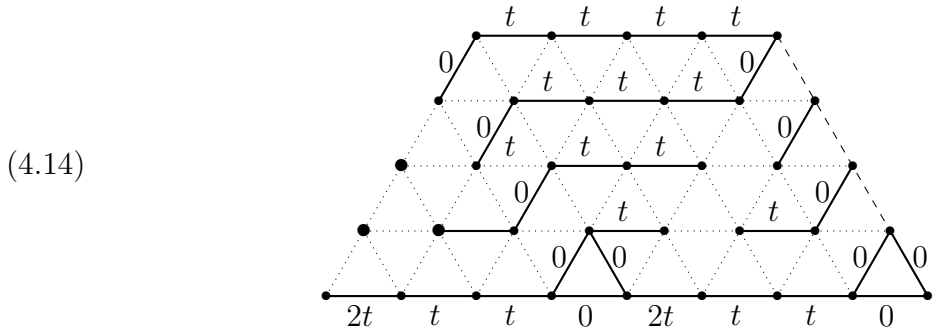
\begin{equation}\label{eqn-Kskel-ex2} 
\vcenter{\hbox{\begin{tikzpicture}[x={(1cm*0.5,-1.732 cm*0.5)},y={(1cm*0.5,1.732 cm*0.5)}] 
\foreach \i in{0,...,4} \draw[dotted](\i,\i)--(\i,4+\i); \foreach \i in{5,...,8} \draw[dotted](\i,\i)--(\i,8);  
\foreach \i in {0,...,3} \draw[thick](\i,3)--(\i,4); \foreach \i in {4,...,7} \draw[thick](\i,7)--(\i,8);
\foreach \j in{1,...,4} \draw[dotted](0,\j)--(\j,\j);  \foreach \j in{5,...,7} \draw[dotted] (\j-4,\j)--(\j,\j); \draw[dashed](4,8)--(8,8); 
\draw[thick](3,4)--(4,4); \draw[thick](7,8)--(8,8); 
\foreach \i in{0,...,4} \draw[dotted](0,\i)--(8-\i,8); 
\draw[thick](1,2)--(2,3); \foreach \i in {1,...,3} \draw[thick](\i,4)--(4,8-\i); \draw[thick](5,6)--(6,7); 
\draw[thick](0,0)--(8,8); \draw[thick](0,4)--(4,8); 
\foreach\j in{0,...,4} \draw(0,\j)node(a0\j){$\sc\bullet$};
\foreach\j in{1,...,5} \draw(1,\j)node(a1\j){$\sc\bullet$};
\foreach\j in{2,...,6} \draw(2,\j)node(a2\j){$\sc\bullet$};
\foreach\j in{3,...,7} \draw(3,\j)node(a3\j){$\sc\bullet$};
\foreach\j in{4,...,8} \draw(4,\j)node(a4\j){$\sc\bullet$};
\foreach\j in{5,...,8} \draw(5,\j)node(a0\j){$\sc\bullet$};
\foreach\j in{6,...,8} \draw(6,\j)node(a1\j){$\sc\bullet$};
\foreach\j in{7,...,8} \draw(7,\j)node(a2\j){$\sc\bullet$};
\foreach\j in{8,...,8} \draw(8,\j)node(a3\j){$\sc\bullet$};
\draw(0,1)node{$\bullet$};\draw(0,2)node{$\bullet$};\draw(1,2)node{$\bullet$};
\foreach\i in{0,...,3} \draw(\i-0.2,4-0.6)node(\i4){$0$};
\foreach\i in{4,...,7} \draw(\i-0.2,8-0.6)node(\i8){$0$};
\draw(4-0.4,4+0.2)node(88){$0$};
\draw(8-0.4,8+0.2)node(88){$0$};
\draw(1-0.4,1-0.7)node(mu1){$2t$};\draw(2-0.4,2-0.7)node(mu2){$t$};\draw(3-0.4,3-0.7)node(mu3){$t$};\draw(4-0.4,4-0.7)node(mu4){$0$};
\draw(2-0.7,4-0.4)node(24){$t$};
\draw(5-0.4,5-0.7)node(nu1){$2t$};\draw(6-0.4,6-0.7)node(nu2){$t$};\draw(7-0.4,7-0.7)node(nu3){$t$};\draw(8-0.4,8-0.7)node(nu4){$0$};
\draw(6-0.7,7-0.4)node(67){$t$};
\foreach\j in{1,...,4} \draw(\j-0.7,4+\j-0.4)node(la\j){$t$};
\foreach\j in{2,...,4} \draw(\j-0.7,3+\j-0.4)node(la\j){$t$};
\foreach\j in{3,...,4} \draw(\j-0.7,2+\j-0.4)node(la\j){$t$};
\foreach\j in{4,...,4} \draw(\j-0.7,1+\j-0.4)node(la\j){$t$};
\end{tikzpicture}}}
\end{equation} 
Here, there are $5$ components disconnected from the horizontal boundaries, but only $3$ of them have been signified by 
an enlarged bullet point at their left-hand extremity. This is because the skeletal upper bound of $5$ in the degree of the
quasi-polynomial $n_{t(2,1,1),t(2,1,1)}^{t(1,1,1,1)}$ can be lowered to $3$ by noting the requirement that the edge labels 
on the right- and left-hand boundaries coincide. This leads to the expectation that the quasi-polynomial is of degree $3$,
as first established in~\cite{GOY1} and confirmed here in (\ref{eqn-GOY}). 

For $\mu=\nu=\lambda=(2,1,1)$ the skeletal $K$-hive is given by: 
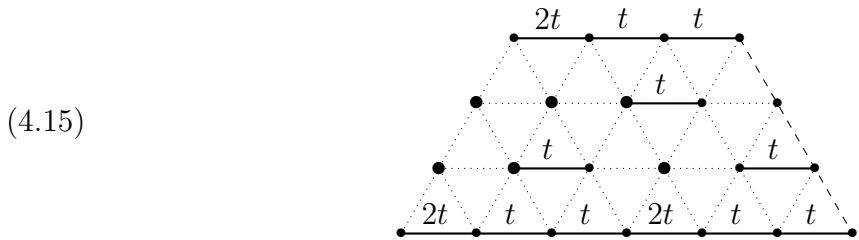
\begin{equation}\label{eqn-NLhive-211-211-211}
\vcenter{\hbox{\begin{tikzpicture}[x={(1cm*0.5,-1.732 cm*0.5)},y={(1cm*0.5,1.732 cm*0.5)}] 
\draw[dotted](0,0)--(0,3); \draw[dotted](1,1)--(1,4); \draw[dotted](2,2)--(2,5); \draw[dotted](3,3)--(3,6); \draw[dotted](4,4)--(4,6);\draw[dotted](5,5)--(5,6);
\draw[dotted](0,1)--(1,1); \draw[dotted](0,2)--(2,2); \draw[dotted](0,3)--(3,3);\draw[dotted](1,4)--(4,4); \draw[dotted](2,5)--(5,5); \draw[dashed](3,6)--(6,6);
\draw[thick](0,0)--(6,6); 
\draw[dotted](0,1)--(5,6);\draw[thick](1,2)--(2,3); \draw[thick](4,5)--(5,6); \draw[thick](2,4)--(3,5); 
\draw[dotted](0,2)--(4,6);
\draw[thick](0,3)--(3,6);
\foreach\j in{0,...,3} \draw(0,\j)node(a0\j){$\sc\bullet$};
\foreach\j in{1,...,4} \draw(1,\j)node(a1\j){$\sc\bullet$};
\foreach\j in{2,...,5} \draw(2,\j)node(a1\j){$\sc\bullet$};
\foreach\j in{3,...,6} \draw(3,\j)node(a2\j){$\sc\bullet$};
\foreach\j in{4,...,6} \draw(4,\j)node(a0\j){$\sc\bullet$};
\foreach\j in{5,...,6} \draw(5,\j)node(a0\j){$\sc\bullet$};
\foreach\j in{6,...,6} \draw(6,\j)node(a0\j){$\sc\bullet$};
\draw(0,1)node(01){$\bullet$};\draw(1,2)node(12){$\bullet$};\draw(3,4)node(34){$\bullet$};\draw(0,2)node(02){$\bullet$};\draw(1,3)node(13){$\bullet$};\draw(2,4)node(24){$\bullet$};
\foreach\j in{1,4} \draw(\j-0.7,\j-0.4)node(mu\j){$2t$};
\foreach\j in{2,3,5,6} \draw(\j-0.7,\j-0.4)node(nu\j){$t$};
\foreach\j in{1}       \draw(\j-0.7,3+\j-0.4)node(la\j){$2t$};
\foreach\j in{2,3}     \draw(\j-0.7,3+\j-0.4)node(la\j){$t$};
\foreach\j in{2,5}     \draw(\j-0.7,1+\j-0.4)node(ka\j){$t$};
\foreach\j in{3}       \draw(\j-0.7,2+\j-0.4)node(ma\j){$t$};
\end{tikzpicture}}}
\end{equation} 
This time there are $8$ components disconnected from the horizontal boundaries. However, the coincidence
of right and left hand boundary edge labels means that there are just $6$ degrees of freedom, leading to the 
expectation that $n_{t(2,1,1)),t(2,1,1)}^{t(2,1,1)}$ is a quasi-polynomial of degree $6$, as confirmed
in (\ref{eqn-GOY}) and first established in a footnote of~\cite{GOY1}. 

\section{Examples and stability phenomena}\label{Sec-examples}

To illustrate cases in which $\mu$ and $\nu$ are fixed and $\lambda$ varies over all partitions
of weight $|\lambda|$ ranging from $0$ to $|\mu|+|\nu|$ it is instructive to consider, for example, 
the case $\mu=\nu=(3,1)$ with $\lambda$ arbitrary. The results are tabulated below, 
stratified by $|\lambda|$ which varies from $8$ down to $0$. One can see that in each case
the positivity of all coefficients in $P_e(t)$, $P_o(t)$ and $G(w)$ is confirmed, thereby supporting 
the parts {\bf E(iv)}, {\bf O(iv)}, {\bf E(v)} and {\bf O(v)} of Conjecture~\ref{Conj-EO}. 


\begin{equation}\label{eqn-3131-part1}  
\begin{array}{|c|c|l|}
\hline
\lambda&N_{(3,1),(3,1)}^\lambda(w)&n_{t(3,1),t(3,1)}^{t\lambda}\cr
\hline\hline
\begin{array}{l}
(62),(61^2),(53),(51^3),(4^2),\cr
(42^2),(421^2),(3^22),(3^21^2)\cr
\end{array}
&\ds \frac{1}{(1-w)}&1\cr
\hline
(521),(431)&\ds \frac{1}{(1-w)^2}&(t+1)\cr
\hline
\hline
(61),(41^3),(32^2),(321^2)&\ds \frac{1}{(1-w^2)}&\begin{array}{cl} 1&t~\mbox{even}\cr
                                             0&t~\mbox{odd} \cr \end{array}
\cr
\hline
(52),(51^2),(43),(3^21)&\ds \frac{1}{(1-w^2)^4}& \begin{array}{cl} (t+2)(t+4)(t+6)/48&t~\mbox{even}\cr
                                             0&t~\mbox{odd} \cr \end{array}
\cr
\hline
(421)&\ds \frac{(1+3w^2)}{(1-w^2)^5}& \begin{array}{cl} (t+2)^2(t+4)(t+6)/96&t~\mbox{even}\cr
                                             0&t~\mbox{odd} \cr \end{array}
\cr
\hline
\hline
(6),(2^3),(2^21^2)&\ds \frac{1}{(1-w)}&1
\cr
\hline
(51),(41^2)&\ds \frac{1}{(1-w)^3(1-w^2)}&\begin{array}{cl} (t+2)(t+4)(2t+3)/24&t~\mbox{even}\cr
                                          (t+1)(t+3)(2t+7)/24&t~\mbox{odd} \cr \end{array}
\cr
\hline
(42)&\ds \frac{1}{(1-w)^4}&(t+1)(t+2)(t+3)/6
\cr
\hline
(3^2)&\ds \frac{1}{(1-w)^2(1-w^2)}&\begin{array}{cl}(t+2)^2/4&t~\mbox{even}\cr
                                           (t+1)(t+3)/4&t~\mbox{odd} \cr \end{array}
\cr
\hline
(321)&\ds \frac{1}{(1-w)^4(1-w^2)}&\begin{array}{cl}(t+2)(t+4)(t^2+6t+6)/48&t~\mbox{even}\cr
                                           (t+1)(t+3)^2(t+5)/48&t~\mbox{odd} \cr \end{array}
\cr
\hline\hline
(5),(21^3)&\ds \frac{1}{(1-w^2)^2}&\begin{array}{cl}1&t~\mbox{even}\cr
                                           0&t~\mbox{odd} \cr \end{array}
\cr
\hline
(41)&\ds \frac{(1+6w^2+4w^4)}{(1-w^2)^4}&\begin{array}{cl}(t+2)(11t^2+26t+24)/48&t~\mbox{even}\cr
                                           0&t~\mbox{odd} \cr \end{array}
\cr
\hline
(32)&\ds \frac{(1+7w^2+5w^4)}{(1-w^2)^4}&\begin{array}{cl}(t+2)(13t^2+28t+24)/48&t~\mbox{even}\cr
                                           0&t~\mbox{odd} \cr \end{array}
\cr
\hline
(31^2)&\ds \frac{(1+3w^2)}{(1-w^2)^4}&\begin{array}{cl}(t+2)(t+4)(2t+3)/24&t~\mbox{even}\cr
                                           0&t~\mbox{odd} \cr \end{array}
\cr
\hline
(2^21)&\ds \frac{(1+w^2)}{(1-w^2)^4}&\begin{array}{cl}(t+2)(t+3)(t+4)/24&t~\mbox{even}\cr
                                           0&t~\mbox{odd} \cr \end{array}
\cr
\hline\hline
\end{array}
\end{equation}

\begin{equation}\label{eqn-3131-part2}  
\begin{array}{|c|c|l|}
\hline
\lambda&N_{(3,1),(3,1)}^\lambda(w)&n_{t(3,1),t(3,1)}^{t\lambda}\cr
\hline\hline
(4)&\ds \frac{1}{(1-w)^2}&(t+1)
\cr
\hline
(31)&\ds \frac{(1+w+w^2)}{(1-w)^3(1-w^2)}&\begin{array}{cl}(t+2)(2t^2+5t+4)/8&t~\mbox{even}\cr
                                           (t+1)(2t^2+7t+7)/8&t~\mbox{odd} \cr \end{array}
\cr
\hline
(2^2)&\ds \frac{1}{(1-w)^3}&(t+1)(t+2)/2
\cr
\hline
(21^2)&\ds \frac{1}{(1-w)^3(1-w^2)}&\begin{array}{cl} (t+2)(t+4)(2t+3)/24&t~\mbox{even}\cr
                                          (t+1)(t+3)(2t+7)/24&t~\mbox{odd} \cr \end{array}
\cr
\hline
(1^4)&\ds \frac{1}{(1-w)}&1
\cr
\hline\hline
(3)&\ds \frac{(1+w^2)}{(1-w^2)^2}&\begin{array}{cl}(t+1)&t~\mbox{even}\cr
                                           0&t~\mbox{odd} \cr \end{array}
\cr
\hline
(21)&\ds \frac{(1+w^2)(1+4w^2)}{(1-w^2)^4}&\begin{array}{cl}(t+2)(5t^2+11t+12)/24&t~\mbox{even}\cr
                                           0&t~\mbox{odd} \cr \end{array}
\cr
\hline
(1^3)&\ds \frac{1}{(1-w^2)^2}&\begin{array}{cl}(t+1)/2&t~\mbox{even}\cr
                                           0&t~\mbox{odd} \cr \end{array}
\cr
\hline
\hline
(2)&\ds \frac{1}{(1-w)^2}&(t+1)
\cr
\hline
(1^2)&\ds \frac{1}{(1-w)^2(1-w^2)}&\begin{array}{cl} (t+2)^2/4&t~\mbox{even}\cr
                                          (t+1)(t+3)/4&t~\mbox{odd} \cr \end{array}
\cr
\hline
\hline
(1)&\ds \frac{1}{(1-w^2)^2}&\begin{array}{cl}(t+2)/2&t~\mbox{even}\cr
                                           0&t~\mbox{odd} \cr \end{array}
\cr
\hline
\hline
(0)&\ds \frac{1}{(1-w)}&1
\cr
\hline
\end{array}
\end{equation}

While continuing to also support the validity of the various positivity conjectures,
an interesting new stability phenomenon shows itself by considering the
effect of adding a fixed integer $a$ to the first part of all three partitions $\mu$, $\nu$ and $\lambda$.
In the case $\mu=(3,3)$, $\nu=(2,1)$ and $\lambda=(3,2)$ one finds
\begin{equation}\label{eqn-stab-ex1} 
\begin{array}{|c|c|l|}
\hline
a&N_{(a\+3,3),(a\+2,1)}^{(a\+3,2)}(w)&~~n_{t(a\+3,3),t(a\+2,1)}^{t(a\+3,2)}\cr
\hline\hline
a=0&\ds \frac{1}{(1\-w)^2(1\-w^2)}&\begin{array}{ll} (t\+2)^2/4&~~t~\mbox{even}\cr
                                         (t\+1)(t\+3)/4&~~t~\mbox{odd}\cr \end{array}
\cr
\hline
a=1&\ds \frac{(6w^4\+8w^2\+1)}{(1\-w^2)^4}&\begin{array}{cl} (t\+2)(5t^2\+10t\+8)/16&~~t~\mbox{even}\cr 
                                           0&~~t~\mbox{odd}\cr \end{array}
\cr
\hline
a=2&\ds \frac{(3w^2\+3w\+1)}{(1\-w)^3(1\-w^2)}&\begin{array}{ll} (t\+2)(14t^2\+23t\+12)/24&~~t~\mbox{even}\cr 
                                           (t\+1)(14t^2\+37t\+21)/24&~~t~\mbox{odd}\cr \end{array}
\cr
\hline

a=3&\ds \frac{(17w^4\+20w^2\+1)}{(1\-w^2)^4}&\begin{array}{cl} (t\+2)(19t^2\+28t\+12)/24&~~t~\mbox{even}\cr 
                                           0&~~t~\mbox{odd}\cr \end{array}
\cr
\hline
\begin{array}{c}
a\geq4\cr
a~\mbox{even}\cr
\end{array}
&\ds \frac{(4w\+1)}{(1\-w)^4}
&\begin{array}{ll}
            (t\+1)(t\+2)(5t\+3)/6&~~\mbox{all $t$}\cr

	\end{array}
\cr
\hline
\begin{array}{c}
a\geq5\cr
a~\mbox{odd}\cr
\end{array}
&\ds \frac{(17w^4\+22w^2\+1)}{(1\-w^2)^4}
&\begin{array}{cl}
           (t\+1)(t\+2)(5t\+3)/6&~~t~\mbox{even}\cr
                          0&~~t~\mbox{odd}\cr 
 \end{array}
\cr																					
\hline
\hline
\end{array}
\end{equation}
Similarly, for $\mu=(3,3)$, $\nu=(2,1)$ and $\lambda=(2,2)$ one finds
\begin{equation}\label{eqn-stab-ex2}  
\begin{array}{|c|c|l|}
\hline
a&N_{(a\+3,3),(a\+2,1)}^{(a\+2,2)}(w)&~~n_{t(a\+3,3),t(a\+2,1)}^{t(a\+2,2)}\cr
\hline\hline
a=0&\ds \frac{1}{(1-w^2)}&\begin{array}{cl}1&~~t~\mbox{even}\cr 
                                           0&~~t~\mbox{odd}\cr \end{array}
\cr
\hline
a=1&\ds \frac{1}{(1\-w)^3(1\-w^2)}&\begin{array}{ll} (t\+2)(t\+4)(2t\+3)/24&~~t~\mbox{even}\cr
                                         (t\+1)(t\+3)(2t\+7)/24&~~t~\mbox{odd}\cr \end{array}
\cr
\hline
a=2&\ds \frac{(7w^4\+11w^2\+1)}{(1\-w^2)^4}&\begin{array}{cl} (t\+2)(19t^2\+40t\+24)/48&~~t~\mbox{even}\cr 
                                           0&~~t~\mbox{odd}\cr \end{array}
\cr
\hline
a=3&\ds \frac{(3w\+1)}{(1\-w)^4}&\begin{array}{ll} (t\+1)(t\+2)(4t\+3)/6&~~\mbox{for all $t$}\cr \end{array}
\cr
\hline
a=4&\ds \frac{(17w^4\+21w^2\+1)}{(1\-w^2)^4}&\begin{array}{cl} (t\+2)(13t^2\+20t\+8)/16&~~t~\mbox{even}\cr 
                                           0&~~t~\mbox{odd}\cr \end{array}
\cr
\hline
\begin{array}{c}
a\geq5\cr
a~\mbox{odd}\cr
\end{array}
&\ds \frac{(4w\+1)}{(1\-w)^4}
&\begin{array}{ll}
            (t\+1)(t\+2)(5t\+3)/6&~~\mbox{all $t$}\cr

	\end{array}
\cr
\hline
\begin{array}{c}
a\geq6\cr
a~\mbox{even}\cr
\end{array}
&\ds \frac{(17w^4\+22w^2\+1)}{(1\-w^2)^4}
&\begin{array}{cl}
           (t\+1)(t\+2)(5t\+3)/6&~~t~\mbox{even}\cr
                          0&~~t~\mbox{odd}\cr 
 \end{array}
\cr																					
\hline
\hline
\end{array}
\end{equation}

These tables include for $a=2$ the illustrative results given in the Introduction.
However, what is particularly striking is the stability of the stretched 
Newell--Littlewood coefficients for sufficiently large $a$, that is for $a\geq4$ in one case 
and $a\geq5$ in the other. A further example of this stability phenomenon 
is provided in the case $\mu=(a\+3,3,2)$, $\nu=(a\+2,1,1)$ and $\lambda=(a\+3,2,1)$
by the following data based on explicit results for $0\leq a\leq 20$:

\begin{footnotesize}
\begin{equation}\label{eqn-stab1}
\begin{array}{|c|l|}
\hline
&\begin{array}{l}
            ~~~~~G(w)/(1\!-\!w)^{d_1}(1\!-\!w^2)^{d_2}\cr
						\begin{cases} P_e(t) & t~\mbox{even}\cr
                          P_o(t) & t~\mbox{odd}\cr
						\end{cases} \cr
 \end{array}\cr
\hline\hline
a=0
&\begin{array}{l}
          ~~~~~1/(1-w)^4(1-w^2)^4\cr
					\begin{cases} (t\!\!+\!\!2)(t\!\!+\!\!4)(t\!\!+\!\!6)(t\!\!+\!\!8)(t\!\!+\!\!10)(t^2\!\!+\!\!12t\!\!+\!\!21)/80840&\cr
											   (t\!\!+\!\!1)(t\!\!+\!\!3)(t\!\!+\!\!5)(t\!\!+\!\!6)(t\!\!+\!\!7)(t\!\!+\!\!11)/80840&\cr
																						\end{cases}\cr
 \end{array}\cr
\hline
a=1
&\begin{array}{l}
         ~~~~~(33w^{10}\!\!+\!\!378w^8\!\!+\!\!823w^6\!\!+\!\!449w^4\!\!+\!\!53w^2\!\!+\!\!1)/(1\-w^2)^9\cr
					\begin{cases} (t\!\!+\!\!2)(t\!\!+\!\!4)(t\!\!+\!\!6)(1737t^5\!\!+\!\!2268t^4\!\!+\!\!113836t^3\!\!+\!\!277488t^2\!\!+\!\!339584t\!\!+\!\!215040)/10321920&\cr
													~~~~~~~~0&\cr
																						\end{cases}\cr
 \end{array}\cr
\hline
a=2
&\begin{array}{l}
         ~~~~~ (20w^6\!\!+\!\!59w^5\!\!+\!\!105w^4\!\!+\!\!86w^3\!\!+\!\!43w^2\!\!+\!\!9w\!\!+\!\!1)/(1\-w)^5(1\-w^2)^4\cr
					\begin{cases} (t\!\!+\!\!2)(t\!\!+\!\!4)(t\!\!+\!\!6)(323t^5\!\!+\!\!3416t^4\!\!+\!\!13886t^3\!\!+\!\!27892t^2\!\!+\!\!29216t\!\!+\!\!13440)/645120&\cr
												(t\!\!+\!\!1)(t\!\!+\!\!3)(t\!\!+\!\!5)(323t^5\!\!+\!\!4385t^4\!\!+\!\!22196t^3\!\!+\!\!54868t^2\!\!+\!\!68777t\!\!+\!\!37611)/645120&\cr
																						\end{cases}\cr
 \end{array}\cr
\hline
a=3
&\begin{array}{l}
         ~~~~~ (161w^{10}\!\!+\!\!1471w^8\!\!+\!\!2676w^6\!\!+\!\!1208w^4\!\!+\!\!113w^2\!\!+\!\!1)/(1\-w^2)^9\cr
					\begin{cases} (t\!\!+\!\!2)(t\!\!+\!\!4)(t\!\!+\!\!6)(2815t^5\!\!+\!\!30932t^4\!\!+\!\!129972t^3\!\!+\!\!260752t^2\!\!+\!\!252032t\!\!+\!\!107520)/5160960&\cr
													~~~~~~~~0&\cr
																						\end{cases}\cr
 \end{array}\cr
\hline
\begin{array}{c}
a\geq 4\cr 
a~~\mbox{even}\cr
\end{array}
&\begin{array}{l}
          ~~~~~(20w^6\!\!+\!\!61w^5\!\!+\!\!113w^4\!\!+\!\!96w^3\!\!+\!\!50w^2\!\!+\!\!11w\!\!+\!\!1)/(1\-w)^5(1\-w^2)^4\cr
					\begin{cases} (t\!\!+\!\!2)^2(t\!\!+\!\!3)(t\!\!+\!\!4)(t\!\!+\!\!6)(22t^3\!\!+\!\!132t^2\!\!+\!\!229t\!\!+\!\!140)/40320&\cr
												(t\!\!+\!\!1)(t\!\!+\!\!2)(t\!\!+\!\!3)(t\!\!+\!\!5)(22t^4\!\!+\!\!264t^3\!\!+\!\!1087t^2\!\!+\!\!1910t\!\!+\!\!1197)/40320&\cr
																						\end{cases}\cr
 \end{array}\cr
\hline
\begin{array}{c}
a\geq 5\cr 
a~~\mbox{odd}\cr
\end{array}
&\begin{array}{l}
          ~~~~~(161w^{10}\!\!+\!\!1471w^8\!\!+\!\!2676w^6\!\!+\!\!1208w^4\!\!+\!\!115w^2\!\!+\!\!1)/(1\-w^2)^9\cr
					\begin{cases} (t\!\!+\!\!2)^2(t\!\!+\!\!3)(t\!\!+\!\!4)(t\!\!+\!\!6)(22t^3\!\!+\!\!132t^2\!\!+\!\!229t\!\!+\!\!140)/40320&\cr
													~~~~~~~~0&\cr
																						\end{cases}\cr
 \end{array}\cr
\hline\hline
\end{array}
\end{equation}
\end{footnotesize}

On the basis of this and many similar results, we offer the following
\begin{Conjecture}\label{Conj-stability-aplus}
For any partitions $\mu$, $\nu$ and $\lambda$, let $\mu=(\mu_1,\sigma)$, $\nu=(\nu_1,\tau)$ and $\lambda=(\lambda_1,\rho)$. 
Then for integer $a$ there exists $N\in\N$ such that $n_{(a+\mu_1,\sigma),(a+\nu_1,\tau)}^{(a+\lambda_1,\rho)}$, 
and more generally $n_{t(a+\mu_1,\sigma),t(a+\nu_1,\tau)}^{t(a+\lambda_1,\rho)}$ for all $t\in\N$,
are independent of $a$ for all $a>N$.
\end{Conjecture}

A similar stability conjecture of the same type is supported by a good deal of evidence of the same type, namely
\begin{Conjecture}\label{Conj-stability-a}
For integer $a$ and any partitions $\sigma$, $\tau$ and $\rho$, 
let $\mu=(a,\sigma)$, $\nu=(a,\tau)$ and $\lambda=(a,\rho)$ with $a\geq\max\{\sigma_1,\tau_1,\rho_1\}$. 
Then there exists $N\in\N$ such that $n_{(a,\sigma),(a,\tau)}^{(a,\rho)}$, and more generally 
$n_{t(a,\sigma),t(a,\tau)}^{t(a,\rho)}$ for all $t\in\N$,
are independent of $a$ for all $a>N$.
\end{Conjecture}

\section{Newell--Littlewood cubes}\label{Sec-cubes}
Since the case $\mu=\nu=\lambda$ appears to be of particular interest~\cite{H}, we gather together here some 
some results and observations about Newell--Littlewood stretched cubes.

For example, $\mu=\nu=\lambda=(a,b)$ with $b=1,2,3$ we find the following data:
\begin{equation}
\begin{array}{|c|c|l|}
\hline
\lambda&N_{\lambda,\lambda}^{\lambda}(w)&\begin{array}{cl} P_e(t) & t~\mbox{even}\cr P_o(t) & t~\mbox{odd}\cr \end{array}\cr
\hline\hline
(1,1)&\ds \frac{1}{(1\-w)(1\-w^2)}&\begin{array}{cl} (t\+2)/2&t~\mbox{even}\cr
                                         (t\+1)/2&t~\mbox{odd}\cr \end{array}
\cr
\hline
(2,1)&\ds \frac{(1\+2w^2)^2}{(1\-w^2)^4}&\begin{array}{cl}(t\+2)(3t^2\+6t\+8)/16&t~\mbox{even}\cr
                                           0&t~\mbox{odd}\cr \end{array}
\cr
\hline
\begin{array}{c}
(a,1)\cr
a\geq3\cr
a~\mbox{odd}\cr
\end{array}
&\ds \frac{(1\+w\+w^2)}{(1\-w)^3(1\-w^2)}
&\begin{array}{cl}
            (t\+2)(2t^2\+5t\+4)/8&t~\mbox{even}\cr
            (t\+1)(2t^2\+7t\+7)/8&t~\mbox{odd}\cr 
	\end{array}
\cr
\hline
\begin{array}{c}
(a,1)\cr 
a\geq4\cr
a~\mbox{even}\cr
\end{array}
&\ds \frac{(1\+7w^2\+4w^4)}{(1\-w^2)^4}
&\begin{array}{cl}
           (t\+2)(2t^2\+5t\+4)/8&t~\mbox{even}\cr
                          0&t~\mbox{odd}\cr 
 \end{array}
\cr																					
\hline
\hline
\end{array}
\end{equation}

\begin{equation}
\begin{array}{|c|c|l|}
\hline
\lambda&N_{\lambda,\lambda}^{\lambda}(w)&\begin{array}{cl} P_e(t) & t~\mbox{even}\cr P_o(t) & t~\mbox{odd}\cr \end{array}\cr
\hline
\hline
(2,2)&\ds \frac{1}{(1\-w)^2}&(t\+1)
\cr
\hline
(3,2)&\ds \frac{(1\+11w\+14w^2\+w^6)}{(1\-w^2)^4}&\begin{array}{cl} (3t\+4)(3t^2\+4t\+4)/16&t~\mbox{even}\cr
                             0&t~\mbox{odd}\cr \end{array}
\cr
\hline
(4,2)&\ds \frac{(1\+2w)^2}{(1\-w)^4}&(t\+1)(3t^2\+3t\+2)/2
\cr
\hline
(5,2)&\ds \frac{(1\+35w^2\+53w^4\+4w^6)}{(1\-w^2)^4}&\begin{array}{cl} (31t^3\+66t^2\+48t\+16)/16&t~\mbox{even}\cr
                             0&t~\mbox{odd}\cr \end{array}
\cr
\hline
\begin{array}{c}
(a,2)\cr 
a\geq6\cr
a~\mbox{even}\cr
\end{array}
&\ds \frac{(1\+7w\+4w^2)}{(1\-w)^4}&(t\+1)(4t^2\+5t\+2)/2\cr																					
\hline
\begin{array}{c}
(a,2)\cr
a\geq7\cr
a~\mbox{odd}\cr
\end{array}
&\ds \frac{(1\+38w^2\+53w^4\+4w^6)}{(1\-w^2)^4}&\begin{array}{cl}
                              (t\+1)(4t^2\+5t\+2)/2&t~\mbox{even}\cr
                             0&t~\mbox{odd}\cr \end{array}
\cr
\hline
\hline
\end{array}
\end{equation}

\begin{equation}
\begin{array}{|c|c|l|}
\hline
\lambda&N_{\lambda,\lambda}^{\lambda}(w)&\begin{array}{cl} P_e(t) & t~\mbox{even}\cr P_o(t) & t~\mbox{odd}\cr \end{array}\cr
\hline
\hline
(3,3)&\ds \frac{(1\+w\+w^2)}{(1\-w)(1\-w^2)}&\begin{array}{cl} (3t\+2)/2&t~\mbox{even}\cr
                             (3t\+1)/2&t~\mbox{odd}\cr \end{array}
\cr
\hline
(4,3)&\ds \frac{(1\+18w^2\+24w^4\+2w^6)}{(1\-w^2)^4}&\begin{array}{cl} (15t^3\+36t^2\+36t\+16)/16&t~\mbox{even}\cr
                                                    0 &t~\mbox{odd}\cr \end{array}
\cr
\hline
(5,3)&\ds \frac{(1\+4w\+w^2)^2}{(1\-w)^3(1\-w^2)}&\begin{array}{cl} (6t^3\+9t^2\+6t\+2)/2&t~\mbox{even}\cr
                                                               (6t^3\+9t^2\+6t\+1)/2&t~\mbox{odd}\cr \end{array}
\cr
\hline
(6,3)&\ds \frac{(1\+72w^2\+150w^4\+20w^6)}{(1\-w^2)^4}&\begin{array}{cl} (3t\+2)(27t^2\+18t\+8)/16&t~\mbox{even}\cr
                                                    0 &t~\mbox{odd}\cr \end{array}
\cr
\hline
(7,3)&\ds \frac{(1\+17w\+36w^2\+20w^3\+w^4)}{(1\-w)^3(1\-w^2)}&\begin{array}{cl} (50t^3\+69t^2\+34t\+8)/8&t~\mbox{even}\cr
                                                                        (50t^3\+69t^2\+34t\+7)/8&t~\mbox{odd}\cr \end{array}
\cr
\hline
(8,3)&\ds \frac{(1\+99w^2\+198w^4\+23w^6)}{(1\-w^2)^4}&\begin{array}{cl} (107t^3\+156t^2\+76t\+16)/16&t~\mbox{even}\cr
                                                    0 &t~\mbox{odd}\cr \end{array}
\cr
\hline
\begin{array}{c}
(a,3)\cr 
a\geq9\cr
a~\mbox{odd}\cr
\end{array}
&\ds \frac{(1\+20w\+39w^2\+20w^3\+w^4)}{(1\-w)^3(1\-w^2)}&\begin{array}{cl} (3t\+2)(18t^2\+15t\+4)/8&t~\mbox{even}\cr
                                                                      (3t\+1)(18t^2\+21t\+7)/8&t~\mbox{odd}\cr \end{array}
\cr
\hline
\begin{array}{c}
(a,3)\cr 
a\geq10\cr
a~\mbox{even}\cr
\end{array}
&\ds \frac{(1\+102w^2\+198w^4\+23w^6)}{(1\-w^2)^4}&\begin{array}{cl} (3t\+2)(18t^2\+15t\+4)/8&t~\mbox{even}\cr
                                                    0 &t~\mbox{odd}\cr \end{array}
\cr
\hline
\hline
\end{array}
\end{equation}

In these three tables, apart from the positivity of all coefficients in $G(w)$, $P_e(t)$ and $P_o(t)$,
the other notable feature is that as a function of $a$ the results for $\mu=\nu=\lambda=(a,b)$
are stable, that is to say independent of $a$, for all even and all odd $a\geq 3b>0$. 
This is borne out by further examples and these stable values may be summarised as follows.
 
\begin{equation}
\begin{array}{|c|c|l|}
\hline
(a,b)~\mbox{with}~a\geq3b>0&N_{(a,b),(a,b)}^{(a,b)}(w)&n_{t(a,b),t(a,b)}^{t(a,b)}\cr
\hline
\hline
a~\mbox{even}~~b~\mbox{even}
&\ds \frac{G_{ee}(w)}{(1\-w)^4}& P_e(t)~~\mbox{for all $t$}\cr
\hline
\begin{array}{c}
a~\mbox{even}~~b~\mbox{odd}\cr
a~\mbox{odd}~~b~\mbox{even}\cr
\end{array}
&\ds \frac{G_{eo}(w)}{(1\-w^2)^4}&\begin{array}{cl} P_e(t) & t~\mbox{even}\cr 0& t~\mbox{odd}\cr \end{array}  \cr
\hline
a~\mbox{odd}~~b~\mbox{odd}
&\ds \frac{G_{oo}(w)}{(1\-w)^3(1\-w^2)}&\begin{array}{cl} P_e(t) & t~\mbox{even}\cr P_o(t) & t~\mbox{odd}\cr \end{array} \cr
\hline\hline
\end{array}
\end{equation}
where 
\begin{equation}\label{eqn-PePo-bt}
    \begin{array}{rcl} 
		      P_e(t)&=&\ds \frac{1}{8} (b\,t+2)(2b^2\,t^2+5b\,t+4)\,,\cr
					\cr
					P_o(t)&=&\ds \frac{1}{8} (b\,t+1)(2b^2\,t^2+7b\,t+7)\,,\cr 
		\end{array}
\end{equation}
and
\begin{equation}
\begin{array}{l}
G_{ee}(w)=1\+(b^3/4\+9b^2/8\+7b/4\-3)\,w\+(b^3\-7b/2\+3)\,w^2\+(b^3/4\-9b^2/8\+7b/4\-1)\,w^3\cr\cr
G_{eo}(w)=1\+(2b^3\+9b^2/2\+7b/2\-3)\,w^2\+(8b^3\-7b\+3)\,w^4\+(2b^3\-9b^2/2\+7b/2\-1)\,w^6 \cr\cr
G_{oo}(w)=1\+(b^3/4\+9b^2/8\+7b/4\-17/8)\,w\+(5b^3/4\+9b^2/8\-7b/4\+3/8)\,w^2 \cr\cr
~~~~~~~~~~~~~~~~~~~~\+(5b^3/4\-9b^2/8\-7b/4\+13/8)\,w^3\+(b^3/4\-9b^2/8\+7b/4\-7/8)\,w^4\,.\cr
\end{array}
\end{equation}
These formulae have been verified for all $a\geq 3b>0$ with $b\leq 9$ and $a\leq 31$.
Quite why the results should be stable for all $a\geq 3b>0$ is not clear, but it may be checked
for each $b>0$ that the stable value is not reached in the case $a=3b-1$.

Moving to length $3$ partitions, in the simplest case we have
\begin{equation}
\begin{array}{|c|c|l|}
\hline
(a,a,a)~~a>0&N_{(a,a,a),(a,a,a)}^{(a,a,a)}(w)&n_{t(a,a,a),t(a,a,a)}^{t(a,a,a)}\cr
\hline
\hline
a~\mbox{even}&\ds \frac{((a/2-1)w+1)}{(1\-w)^2}& (at+2)/2~~\mbox{for all $t$}\cr
\hline
a~\mbox{odd}&\ds \frac{((a-1)w^2+1)}{(1-w^2)^2}&\begin{array}{cl} (at+2)/2& t~\mbox{even}\cr 0& t~\mbox{odd}\cr \end{array}  \cr
\hline
\hline
\end{array}
\end{equation}
As a less trivial example we offer:
\begin{equation}
\begin{array}{|c|c|l|}
\hline
(a,b,b)~\mbox{with}~a\geq3b>0&N_{(a,b,b),(a,b,b)}^{(a,b,b)}(w)&n_{t(a,b,b),t(a,b,b)}^{t(a,b,b)}\cr
\hline
\hline
a~\mbox{even}~~b~\mbox{even}
&\ds \frac{G_{ee}(w)}{(1\-w)^7}& P_e(t)~~\mbox{for all $t$}\cr
\hline
a~\mbox{even}~~b~\mbox{odd}
&\ds \frac{G_{eo}(w)}{(1\-w)^3(1\-w^2)^4}&\begin{array}{cl} P_e(t) & t~\mbox{even}\cr P_o(t)& t~\mbox{odd}\cr \end{array}  \cr
\hline
a~\mbox{odd}~~b~\mbox{even}
&\ds \frac{G_{oe}(w)}{(1\-w^2)^7}&\begin{array}{cl} P_e(t) & t~\mbox{even}\cr 0 & t~\mbox{odd}\cr \end{array} \cr
\hline
a~\mbox{odd}~~b~\mbox{odd}
&\ds \frac{G_{oo}(w)}{(1\-w^2)^7}&\begin{array}{cl} P_e(t) & t~\mbox{even}\cr P_o(t) & t~\mbox{odd}\cr \end{array} \cr
\hline\hline
\end{array}
\end{equation}
where 
\begin{equation}\label{eqn-PePo-bt6}
    \begin{array}{rcl} 
		      P_e(t)&=&\ds \frac{1}{1920} (bt+2)(bt+4)(4b^4t^4+36b^3t^3+137b^2t^2+270bt+240)\,,\cr
					\cr
					P_o(t)&=&\ds \frac{1}{1920} (bt+1)(bt+3)(4b^4t^4+44b^3t^3+197b^2t^2+400bt+315)\,,\cr 
		\end{array}
\end{equation}
with 
\begin{equation}
\begin{array}{|c|l|l|}
\hline
G(w)&\mbox{Degree of $G(w)$}&(d_1,d_2)\cr
\hline
G_{ee}(w)&6&(7,0) \cr
G_{eo}(w)&10&(3,4) \cr
G_{oe}(w)&12&(0,7)\cr
G_{oo}(w)&12&(0,7)\cr
\hline
\end{array}
\end{equation}
where $G_{ee}(w)$ and $G_{eo}(w)$ are polynomials in $w$, while $G_{oe}(w)$ and $G_{oo}(w)$ are polynomials in $w^2$, all with positive coefficients.
Again, quite why the results should be stable for all $a\geq 3b>0$ is not clear, but it may be verified 
for each $b>0$ that the stable value is not reached in the case $a=3b-1$.

We can also offer a case for which we have a two-parameter explicit formula: 
\begin{equation}
\begin{array}{|c|c|l|}
\hline
(a,a,b)~\mbox{with}~a\geq3b>0&N_{(a,a,b),(a,a,b)}^{(a,a,b)}(w)&n_{t(a,a,b),t(a,a,b)}^{t(a,a,b)}\cr
\hline
\hline
a~\mbox{even}~~b~\mbox{even}
&\ds \frac{G_{ee}(w)}{(1\-w)^7}& P_e(t)~~\mbox{for all $t$}\cr
\hline
a~\mbox{even}~~b~\mbox{odd}
&\ds \frac{G_{eo}(w)}{(1\-w^2)^7}&\begin{array}{cl} P_e(t) & t~\mbox{even}\cr 0& t~\mbox{odd}\cr \end{array}  \cr
\hline
a~\mbox{odd}~~b~\mbox{even}
&\ds \frac{G_{oe}(w)}{(1\-w)^3(1\-w^2)^4}&\begin{array}{cl} P_e(t) & t~\mbox{even}\cr P_o(t) & t~\mbox{odd}\cr \end{array} \cr
\hline
a~\mbox{odd}~~b~\mbox{odd}
&\ds \frac{G_{oo}(w)}{(1\-w^2)^7}&\begin{array}{cl} P_e(t) & t~\mbox{even}\cr 0 & t~\mbox{odd}\cr \end{array} \cr
\hline\hline
\end{array}
\end{equation}
where 
\begin{equation}\label{eqn-PePo-abt}
    \begin{array}{rcl} 
		      P_e(t)&=&\ds \frac{1}{7680} (bt\+2)((108a\-179)b^5t^5\+(864a/b\-1412)b^4t^4\+(2592a/b\-3916)b^3t^3\cr
					      &&~~~~~~~~~~~~~~~~~+(3456a/b\-3808)b^2t^2+(1920a/b\+960)bt+3840)\,,\cr
					\cr
					P_o(t)&=&\ds \frac{1}{7680} (bt\+2)((108a\-179)b^5t^5\+(864a/b\-1412)b^4t^4\+(2592a/b\-3916)b^3t^3\cr
					      &&~~~~~~~~~~~~~~~~~+(3456a/b\-3568)b^2t^2+(1920a/b)bt+1920)\,,\cr
		\end{array}
\end{equation}
with $G_{ee}(w)$ a polynomial in $w$ of degree $6$, $G_{eo}(w)$ and $G_{oo}(w)$ polynomials in $w^2$ of degree $5$, and $G_{oe}(w)$ a 
polynomial in $w$ of degree $9$, all with positive coefficients. It can be checked that because $a\geq3b\geq0$, all coefficients
in $P_e(t)$ and $P_o(t)$ are positive. These polynomials are of degree $6$, as would be predicted by means of the hive model.

In the most general length $3$ cubic case for which $\mu=\nu=\lambda=(a,b,c)$ we do not have any stable limit $3$-parameter formula, 
but we gather together some data in Appendix~\ref{Sec-appendix} that give a striking demonstration of the positivity of 
coefficients in $P_e(t)$, $P_o(t)$ anf $G(w)$.

\section{Conclusion}\label{Sec-conclusion}
By noting that the Newell--Littlewood coefficients $n_{\mu,\nu}^{\lambda}$ are nothing other than the Clebsch--Gordon coefficients $m_{\mu,\nu}^{\lambda}(\g)$
for orthogonal and symplectic Lie algebras $\g$ of sufficiently high rank, and leaning on the work of De Loera and McAllister~\cite{DeLMcA} it has been established 
in Proposition~\ref{Prop-quasipol} that stretched Newell--Littlewood coefficients are quasi-polynomial of minimum quasi-period at most $2$. 
Taking into account the evenness or oddness of $|\mu|+|\nu|+\|\lambda|$ this has as corollary the validity of two of the items, {\bf E(iii)} and {\bf O(iii)},
included in the list of items in Conjecture~\ref{Conj-EO}. In a somewhat similar manner, the prior work of Belkale and Kumar~\cite{BK} and Sam~\cite{Sa} has 
led to a proof of the saturation property {\bf O(i)}, but not {\bf E(i)} for which there is still room for a counterexample to arise. 

Unfortunately, attempts to provide proofs of the validity of these and the other remaining items of Conjecture~\ref{Conj-EO}, 
are beyond the scope of this particular study. However, thanks to the 
use of universal characters, a constant term formula has been established in Theorem~\ref{The-GF} for the generating function, $N_{\mu,\nu}^\lambda(w)$,
of the stretched Newell--Littlewood coefficients $n_{t\mu,t\nu}^{t\lambda}$. This has enabled both the generating function itself and the 
corresponding quasi-polynomials to be evaluated explicitly without the need for any fitting of data. Beyond the confirmation of results
both established and conjectured in~\cite{DeLMcA}, this has allowed the accumulation of a good deal of evidence in support of the positivity conjectures 
forming parts {\bf E(iv)}, {\bf O(iv)}, {\bf E(v)} and {\bf O(v)} of Conjecture~\ref{Conj-EO}. 

The hive model of Section~\ref{Sec-hives} offers an alternative method of calculating Newell--Littlewood coefficients using Proposition~\ref{Prop-NLcoeff-hive}.
This model has led to a proof in Theorem~\ref{The-Ehrhart} that stretched Newell--Littlewood coefficients are the Erhart quasi-polynomials associated with a 
integer convex polytope specified by means of triangular and rhombus conditions, (\ref{eqn-triangle-conditions}) and (\ref{eqn-rhombus-conditions}). This polytope 
is much simpler to specify than that used for stretched Clebsch--Gordon coefficients~\cite{BZ2,DeLMcA}. In fact it is much more akin to the Newell--Littlewood
polytope introduced by Gao {\it et al.}~\cite{GOY1}, differing only in its emphasis on vertex rather than edge labels in the underlying hive.  
No attempt has been made to establish the precise 
conditions on $\mu$, $\nu$ and $\lambda$ under which $n_{\mu,\nu}^{\lambda}$ is non-vanishing, but this has been pursued by Gao~{\it et al.}~\cite{GOY2}. 

The accumulated data has shown the existence of somewhat unexpected stability phenomena that form the basis of Conjectures~\ref{Conj-stability-aplus} and
\ref{Conj-stability-a}. The latter is reminiscent of the stability property Kronecker of coefficients that led to the concept of reduced Kronecker
coefficients~\cite{M,K3}. This can be seen by setting $|\mu|=|\nu|=|\lambda|=m$ in Conjecture~\ref{Conj-stability-aplus}. This then
yields as a special case 
\begin{Conjecture}\label{Conj-reduced}
For any partitions $\sigma$, $\tau$ and $\rho$, and positive integer $m$, there exists $N\in\N$ such that both
$n_{(m-|\sigma|,\sigma),(m-|\tau|,\tau)}^{(m-|\rho|,\rho)}$ and $n_{t(m-|\sigma|,\sigma),t(m-|\tau|,\tau)}^{t(m-|\rho|,\rho)}$ are independent of $m$ for all $m>N$.
\end{Conjecture}
The generating function approach by means of Theorem~\ref{The-GF} offers the possibility of firmly establishing the notion of such
$m$-independent reduced coefficients in the Newell--Littlewood context. The data gathered together on stretched Newell--Littlewood 
cubes in Section~\ref{Sec-cubes} and the Appendix~\ref{Sec-appendix} points to the possibility of being able to fix the point at which stability 
is reached, at least in this case of cubes.

While no direct attack has been made on the analogs of the Fulton conjecture~\cite{KTW} as offered in 
parts {\bf E(ii)} and {\bf O(ii)} of Conjecture~\ref{Conj-EO}, the data in support of the various positivity conjectures also
supports these. In fact, along with {\bf E(i)} and {\bf O(i)}, they should be rather easy corollaries to the positivity 
conjectures. It is intended, and indeed hoped that both the methods described here and the accumulated data will serve as an 
impetus to further work on these problems.

\begin{landscape}

\begin{appendix}
\section{Appendix}\label{Sec-appendix}

The following data on $n_{t\mu,t\nu}^{t\lambda}$ in the case $\mu=\nu=\lambda=(a,b,c)$ with $a>b>c>0$
supports the validity of parts {\bf E(iv)} and {\bf O(iv)} of Conjecture~\ref{Conj-EO} concerning the positivity of 
the quasi-polynomial coefficients, as well as the stability Conjecture~\ref{Conj-stability-a} arising for all $a\geq3b$.
This data has been gathered by evaluating $N_{\mu,\nu}^{\lambda}(w)$ in the form $G(w)/(1-w)^{d_1}(1-w^2)^{d_2}$
through the use of Theorem~\ref{The-GF},
and in so doing checking that in every instance $G(w)$ is a polynomial with positive integer coefficients in accordance with 
{\bf E(v)} and {\bf O(v)} of Conjecture~\ref{Conj-EO}.

In the case $(a,2,1)$ of the simplest cubes of length $3$ with $a>b>c>0$ we find:

\begin{footnotesize}
\begin{equation}\nonumber
\begin{array}{|c|l|}
\hline
(a,b,c)&\begin{array}{ll} P_e(t) & t~\mbox{even}\cr  P_o(t) & t~\mbox{odd}\cr \end{array}\cr
\hline\hline
(3,2,1)&
\begin{array}{l}
        (t\!\!+\!\!2)(448t^8\!\!+\!\!5401t^7\!\!+\!\!28972t^6\!\!+\!\!91870t^5\!\!+\!\!191110t^4\!\!+\!\!272728t^3\!\!+\!\!272760t^2\!\!+\!\!183456t\!\!+\!\!80640)/161280\cr
			 (t\!\!+\!\!3)(t\!\!+\!\!1)(448t^7\!\!+\!\!4505t^6\!\!+\!\!20410t^5\!\!+\!\!54659t^4\!\!+\!\!94984t^3\!\!+\!\!111035t^2\!\!+\!\!83454t\!\!+\!\!33705)/161280\cr
					\end{array}\cr
\hline
(4,2,1)&
\begin{array}{l} 
					 (t\!\!+\!\!2)(55987t^8\!\!+\!\!703627t^7\!\!+\!\!3826570t^6\!\!+\!\!11870644t^5\!\!+\!\!23340040t^4\!\!+\!\!30448384t^3\!\!+\!\!26725248t^2\!\!+\!\!15344640t\!\!+\!\!5160960)/10321920\cr
					 ~~~~~~~0\cr
\end{array}\cr
\hline
(5,2,1)&
\begin{array}{l}
         
					(t\!\!+\!\!2)(5441t^8\!\!+\!\!69758t^7\!\!+\!\!389318t^6\!\!+\!\!1243412t^5\!\!+\!\!2507732t^4\!\!+\!\!3313136t^3\!\!+\!\!2898864t^2\!\!+\!\!1635264t\!\!+\!\!483840)/967680\cr
					(t\!\!+\!\!3)(t\!\!+\!\!1)(5441t^7\!\!+\!\!58876t^6\!\!+\!\!277007t^5\!\!+\!\!737392t^4\!\!+\!\!1213967t^3\!\!+\!\!1260556t^2\!\!+\!\!793761t\!\!+\!\!249480)/967680\cr
\end{array}\cr
\hline
\begin{array}{c}(a,2,1)\cr a\geq6~\mbox{even}\cr\end{array}
&\begin{array}{l}
					(t\!\!+\!\!2)(10883t^8\!\!+\!\!139559t^7\!\!+\!\!779384t^6\!\!+\!\!2493518t^5\!\!+\!\!5047520t^4\!\!+\!\!6700760t^3\!\!+\!\!5845248t^2\!\!+\!\!3212928t\!\!+\!\!967680)/1935360\cr
					 ~~~~~~~0\cr
					\end{array}\cr
\hline
\begin{array}{c}(a,2,1)\cr a\geq7~\mbox{odd}\cr\end{array}
&\begin{array}{l}
					(t\!\!+\!\!2)(10883t^8\!\!+\!\!139559t^7\!\!+\!\!779384t^6\!\!+\!\!2493518t^5\!\!+\!\!5047520t^4\!\!+\!\!6700760t^3\!\!+\!\!5845248t^2\!\!+\!\!3212928t\!\!+\!\!967680)/1935360\cr
					(t\!\!+\!\!3)(t\!\!+\!\!1)(10883t^7\!\!+\!\!117793t^6\!\!+\!\!554681t^5\!\!+\!\!1480183t^4\!\!+\!\!2449781t^3\!\!+\!\!2556127t^2\!\!+\!\!1576527t\!\!+\!\!446985)/1935360\cr
\end{array}\cr
\hline\hline
\end{array}
\end{equation}
\end{footnotesize}

Here, and in all other cases pursued to date, we find that a stable limit is reached for all $a\geq3b$ with $b>c>0$. 
In such cases we find:

\begin{footnotesize}
\begin{equation}\nonumber
\begin{array}{|c|c|c|}
\hline
(a,b,c)~\mbox{with}~a\geq 3b&N_{(a,b,c),(a,b,c)}^{(a,b,c)}(w)&n_{t(a,b,c),t(a,b,c)}^{t(a,b,c)}\cr
\hline
\hline
\begin{array}{c}a+b+c~\mbox{even}\cr b,c~\mbox{both even}\cr\end{array}
&\ds \frac{G_{ee}(w)}{(1\-w)^{10}}&P_e(t)~~\mbox{for all $t$}\cr
\hline
\begin{array}{c}a+b+c~\mbox{even}\cr b,c~\mbox{not both even}\cr\end{array}
&\ds \frac{G_{eo}(w)}{(1\-w)^6(1\-w^2)^4}&\begin{array}{cl} P_e(t) & t~\mbox{even}\cr P_o(t) & t~\mbox{odd}\cr \end{array} \cr
\hline
a+b+c~\mbox{odd}
&\ds \frac{G_{oo}(w)}{(1\-w^2)^{10}}&\begin{array}{cl} P_e(t) & t~\mbox{even}\cr 0& t~\mbox{odd}\cr \end{array}  \cr
\hline
\hline\hline
\end{array}
\end{equation}
\end{footnotesize}
with
\begin{footnotesize}
\begin{equation}\nonumber
\begin{array}{|c|l|}
\hline
(a,b,c)&~~\mbox{Either}~~\begin{array}{l} P_e(t)\cr P_o(t)\cr  \end{array}~~~~~\mbox{or}~~~P_e(t)~~\mbox{if $P_e(t)=P_o(t)$}\cr
\hline\hline
(a,2,1)&\begin{array}{l}
					(t\!\!+\!\!2)(10883t^8\!\!+\!\!139559t^7\!\!+\!\!779384t^6\!\!+\!\!2493518t^5\!\!+\!\!5047520t^4\!\!+\!\!6700760t^3\!\!+\!\!5845248t^2\!\!+\!\!3212928t\!\!+\!\!967680)/1935360\cr
					(t\!\!+\!\!3)(t\!\!+\!\!1)(10883t^7\!\!+\!\!117793t^6\!\!+\!\!554681t^5\!\!+\!\!1480183t^4\!\!+\!\!2449781t^3\!\!+\!\!2556127t^2\!\!+\!\!1576527t\!\!+\!\!446985)/1935360\cr
					\end{array}\cr
\hline
(a,3,1)&\begin{array}{l}
     (t\!\!+\!\!2)(7699t^8\!\!+\!\!88300t^7\!\!+\!\!421804t^6\!\!+\!\!1101286t^5\!\!+\!\!1742992t^4\!\!+\!\!1741588t^3\!\!+\!\!1095000t^2\!\!+\!\!410976t\!\!+\!\!80640)/161280\cr
     (t\!\!+\!\!3)(t\!\!+\!\!1)(7699t^7\!\!+\!\!72902t^6\!\!+\!\!283699t^5\!\!+\!\!591392t^4\!\!+\!\!728899t^3\!\!+\!\!537800t^2\!\!+\!\!220119t\!\!+\!\!37170)/161280\cr
					\end{array}\cr
\hline
(a,3,2)&\begin{array}{l}
     (t\!\!+\!\!1)(190727t^8\!\!+\!\!1525816t^7\!\!+\!\!5461022t^6\!\!+\!\!11404708t^5\!\!+\!\!15266888t^4\!\!+\!\!13518064t^3\!\!+\!\!7967808t^2\!\!+\!\!3036672t\!\!+\!\!645120)/645120\cr
     (t+1)^3(190727t^6\!\!+\!\!1144362t^5\!\!+\!\!2981571t^4\!\!+\!\!4297204t^3\!\!+\!\!3690909t^2\!\!+\!\!1839042t\!\!+\!\!452025)/645120\cr    
					\end{array}\cr
\hline
(a,4,1)&\begin{array}{l}
     (t\!\!+\!\!2)(26011t^8\!\!+\!\!286576t^7\!\!+\!\!1288150t^6\!\!+\!\!3076144t^5\!\!+\!\!4302184t^4\!\!+\!\!3661564t^3\!\!+\!\!1884120t^2\!\!+\!\!552096t\!\!+\!\!80640)/161280\cr
     (t\!\!+\!\!3)(t\!\!+\!\!1)(26011t^7\!\!+\!\!234554t^6\!\!+\!\!845053t^5\!\!+\!\!1568570t^4\!\!+\!\!1645033t^3\!\!+\!\!980090t^2\!\!+\!\!308319t\!\!+\!\!37170)/161280\cr
					\end{array}\cr
\hline
(a,4,2)&\begin{array}{l}
   (t\!\!+\!\!1)(43532t^8\!\!+\!\!279118t^7\!\!+\!\!779384t^6\!\!+\!\!1246759t^5\!\!+\!\!1261880t^4\!\!+\!\!837595t^3\!\!+\!\!365328t^2\!\!+\!\!100404t\!\!+\!\!15120)/15120\cr
 					\end{array}\cr
\hline
(a,4,3)&\begin{array}{l}
        (1860257t^9\!\!+\!\!13261743t^8\!\!+\!\!42512874t^7\!\!+\!\!80444826t^6\!\!+\!\!99261036t^5\!\!+\!\!83191752t^4\!\!+\!\!48076336t^3\!\!+\!\!19027584t^2\!\!+\!\!4891392t\!\!+\!\!645120)/645120\cr
        (t\!\!+\!\!1)^2(1860257t^7\!\!+\!\!9541229t^6\!\!+\!\!21570159t^5\!\!+\!\!27763279t^4\!\!+\!\!22164319t^3\!\!+\!\!11099835t^2\!\!+\!\!3289617t\!\!+\!\!446985)/645120\cr
\end{array}\cr 
\hline
(a,5,1)&\begin{array}{l}
(t\!\!+\!\!2)(61387t^8\!\!+\!\!661900t^7\!\!+\!\!2876944t^6\!\!+\!\!6523006t^5\!\!+\!\!8442784t^4\!\!+\!\!6438628t^3\!\!+\!\!2854680t^2\!\!+\!\!693216t\!\!+\!\!80640)/161280\cr
(t\!\!+\!\!3)(t\!\!+\!\!1)(61387t^7\!\!+\!\!539126t^6\!\!+\!\!1860079t^5\!\!+\!\!3219200t^4\!\!+\!\!3031759t^3\!\!+\!\!1539560t^2\!\!+\!\!396519t\!\!+\!\!37170)/161280 \cr
					\end{array}\cr
\hline
(a,5,2)&\begin{array}{l}
(19719167t^9\!\!+\!\!137161701t^8\!\!+\!\!413898366t^7\!\!+\!\!714090510t^6\!\!+\!\!781934748t^5\!\!+\!\!568830024t^4\!\!+\!\!278307184t^3\!\!+\!\!90083520t^2\!\!+\!\!18224640t\!\!+\!\!1935360)/1935360\cr
(t\!\!+\!\!1)(19719167t^8\!\!+\!\!117442534t^7\!\!+\!\!296455832t^6\!\!+\!\!417634678t^5\!\!+\!\!364300070t^4\!\!+\!\!204529954t^3\!\!+\!\!73658160t^2\!\!+\!\!16079490t\!\!+\!\!1696275)/1935360\cr
		\end{array}\cr
\hline
(a,5,3)&\begin{array}{l}
(14374271t^9\!\!+\!\!78323328t^8\!\!+\!\!192117936t^7\!\!+\!\!278601120t^6\!\!+\!\!264107886t^5\!\!+\!\!170634492t^4\!\!+\!\!76100344t^3\!\!+\!\!23213520t^2\!\!+\!\!4609728t\!\!+\!\!483840)/483840\cr
(t\!\!+\!\!1)(14374271t^8\!\!+\!\!63949057t^7\!\!+\!\!128168879t^6\!\!+\!\!150432241t^5\!\!+\!\!113675645t^4\!\!+\!\!56958847t^3\!\!+\!\!18294777t^2\!\!+\!\!3308463t\!\!+\!\!238140))/483840\cr
		\end{array}\cr
\hline
(a,5,4)&\begin{array}{l}
(9374843t^9\!\!+\!\!57962079t^8\!\!+\!\!159266238t^7\!\!+\!\!255773490t^6\!\!+\!\!265478220t^5\!\!+\!\!185644536t^4\!\!+\!\!88964272t^3\!\!+\!\!29037120t^2\!\!+\!\!6100992t\!\!+\!\!645120)/645120\cr
(t\!\!+\!\!1)^2(9374843t^7\!\!+\!\!39212393t^6\!\!+\!\!71466609t^5\!\!+\!\!73627879t^4\!\!+\!\!46755853t^3\!\!+\!\!18504951t^2\!\!+\!\!4298247t\!\!+\!\!452025)/645120\cr
		\end{array}\cr
\hline
(a,6,1)&\begin{array}{l}
(t\!\!+\!\!2)(119371t^8\!\!+\!\!1269712t^7\!\!+\!\!5401126t^6\!\!+\!\!11832472t^5\!\!+\!\!14509528t^4\!\!+\!\!10193740t^3\!\!+\!\!4006680t^2\!\!+\!\!834336t\!\!+\!\!80640)/161280\cr
(t\!\!+\!\!3)(t\!\!+\!\!1)(119371t^7\!\!+\!\!1030970t^6\!\!+\!\!3458557t^5\!\!+\!\!5707586t^4\!\!+\!\!4968457t^3\!\!+\!\!2216210t^2\!\!+\!\!484719t\!\!+\!\!37170))/161280\cr
		\end{array}\cr
\hline
(a,6,2)&\begin{array}{l}
(t\!\!+\!\!1)(61592t^8\!\!+\!\!353200t^7\!\!+\!\!843608t^6\!\!+\!\!1101286t^5\!\!+\!\!871496t^4\!\!+\!\!435397t^3\!\!+\!\!136875t^2\!\!+\!\!25686t\!\!+\!\!2520)/2520\cr
		\end{array}\cr
\hline
(a,6,3)&\begin{array}{l}
     (3t\!\!+\!\!2)(2644569t^8\!\!+\!\!11304279t^7\!\!+\!\!21043368t^6\!\!+\!\!22441662t^5\!\!+\!\!15142560t^4\!\!+\!\!6700760t^3\!\!+\!\!1948416t^2\!\!+\!\!356992t\!\!+\!\!35840)/71680\cr
		(3t\!\!+\!\!1)(t\!\!+\!\!1)(2644569t^7\!\!+\!\!9541233t^6\!\!+\!\!14976387t^5\!\!+\!\!13321647t^4\!\!+\!\!7349343t^3\!\!+\!\!2556127t^2\!\!+\!\!525509t\!\!+\!\!49665)/71680
\end{array}\cr
\hline
(a,6,4)&\begin{array}{l}
(2t\!\!+\!\!1)(381454t^8\!\!+\!\!1525816t^7\!\!+\!\!2730511t^6\!\!+\!\!2851177t^5\!\!+\!\!1908361t^4\!\!+\!\!844879t^3\!\!+\!\!248994t^2\!\!+\!\!47448t\!\!+\!\!5040)/5040\cr
		\end{array}\cr
\hline
(a,6,5)&\begin{array}{l}
(33064469t^9\!\!+\!\!185217615t^8\!\!+\!\!455264850t^7\!\!+\!\!646494282t^6\!\!+\!\!587226444t^5\!\!+\!\!356073480t^4\!\!+\!\!146931760t^3\!\!+\!\!41033088t^2\!\!+\!\!7310592t\!\!+\!\!645120)/645120\cr
(t\!\!+\!\!1)^2(33064469t^7\!\!+\!\!119088677t^6\!\!+\!\!184023027t^5\!\!+\!\!159359551t^4\!\!+\!\!84484315t^3\!\!+\!\!27745299t^2\!\!+\!\!5296797t\!\!+\!\!446985)/645120\cr
\end{array}\cr
\hline
\end{array}
\end{equation}
\end{footnotesize}

\end{appendix}

\end{landscape}

\bigskip \bigskip

\noindent{\bf Acknowledgement} \medskip

This work was stimulated by the work of Gao, Orelowitz and Yong~\cite{GOY1}, and it is a pleasure to thank Professor Yong for the
subsequent exchange of correspondence and for his encouragement. Further thanks are due to Professor Sam for bringing to my attention
his saturation theorem for the classical groups in~\cite{Sa} and for additional helpful remarks.

\end{document}